\documentclass[12pt,reqno]{amsart}
\usepackage{enumerate, latexsym, amsmath, amsfonts, amssymb, amsthm, color}
\def\pmod #1{\ ({\rm{mod}}\ #1)}
\def\Z{\Bbb Z}
\def\N{\Bbb N}

\def\Q{\Bbb Q}

\def\R{\Bbb R}

\def\l{\left}
\def\r{\right}
\def\bg{\bigg}
\def\({\bg(}
\def\){\bg)}
\def\t{\text}
\def\f{\frac}
\def\mo{{\rm{mod}\ }}

\def\ls{\leqslant}
\def\gs{\geqslant}

\def\sm{\setminus}
\def\bi{\binom}
\def\al{\alpha}

\def\ve{\varepsilon}

\def\eq{\equiv}

\def\da{\delta}

\def\Proof{\noindent{\it Proof}}
\def\Ack{\medskip\noindent {\bf Acknowledgments}}
\theoremstyle{plain}
\newtheorem{theorem}{Theorem}

\newtheorem{lemma}{Lemma}
\newtheorem{corollary}{Corollary}
\newtheorem{conjecture}{Conjecture}
\theoremstyle{definition}

\theoremstyle{remark}
\newtheorem{remark}{Remark}

 \vspace{4mm}

\begin{document}

\hbox{Included in: {\it Combinatorial and Additive Number Theory: CANT 2011 and 2012}}
\hbox{(edited by M.B. Nathanson), Springer Proc. in Math. \& Stat., Vol. 101,}
\hbox{Springer, New York, 2014, pp. 257--312.}
\medskip

\title
[{Central binomial and trinomial coefficients}]
{On sums related to central binomial and trinomial coefficients}

\author
[Zhi-Wei Sun] {Zhi-Wei Sun}

\address {Department of Mathematics, Nanjing
University, Nanjing 210093, People's Republic of China}
\email{zwsun@nju.edu.cn}


 \begin{abstract}  A generalized central trinomial coefficient $T_n(b,c)$ is
the coefficient of $x^n$ in the expansion of $(x^2+bx+c)^n$ with
$b,c\in\Z$. In this paper we investigate congruences and series for
sums of terms related to  central binomial coefficients and
generalized central trinomial coefficients.
The paper contains many conjectures on congruences related to
representations of primes by certain binary quadratic forms, and 62
proposed new series for $1/\pi$ motivated by congruences and related
dualities.
\end{abstract}

\maketitle

\section{Introduction}
\setcounter{lemma}{0}
\setcounter{theorem}{0}
\setcounter{corollary}{0}
\setcounter{remark}{0}
\setcounter{equation}{0}
\setcounter{conjecture}{0}

Let $\N=\{0,1,2,\ldots\}$. The central binomial coefficients
$$\bi{2n}n=\f{(2n)!}{(n!)^2}\quad (n\in\N)$$
play important roles in combinatorics and number theory.
In this section we first review some known results on sums involving products of at most three central binomial coefficients.

Let $\Z^+=\{1,2,3,\ldots\}$. Recall that for given numbers $A$ and $B$ the Lucas sequence
$u_n=u_n(A,B)\ (n\in\N)$ and its companion $v_n=v_n(A,B)\ (n\in\N)$
are defined by
$$u_0=0,\ u_1=1,\ u_{n+1}=Au_n-Bu_{n-1}\ (n\in\Z^+),$$
and
$$v_0=2,\ v_1=A,\ v_{n+1}=Av_n-Bv_{n-1}\ (n\in\Z^+).$$
It is well known that
$$(\al-\beta)u_n=\al^n-\beta^n\ \ \t{and}\ \ v_n=\al^n+\beta^n\quad\t{for all}\ n\in\N,$$
where $\al=(A+\sqrt{\Delta})/2$ and $\beta=(A-\sqrt{\Delta})/2$ with
$\Delta=A^2-4B$.

Let $p$ be an odd prime and let $m$ be any integer not divisible by $p$.
The author [Su1] proved that
$$\sum_{k=0}^{p-1}\f{\bi{2k}k}{m^k}\eq\l(\f{m^2-4m}p\r)+u_{p-(\f{m^2-4m}p)}(m-2,1)\ (\mo\ p^2),$$
where $(-)$ denotes the Jacobi symbol.

Let $p\eq1\ (\mo\ 4)$ be a prime. Write $p=x^2+y^2$ with $x\eq1\
(\mo\ 4)$ and $y\eq0\ (\mo\ 2)$. Gauss showed in 1828 the famous
congruence
$$\bi{(p-1)/2}{(p-1)/4}\eq2x\pmod{p},$$
and this was further refined by S. Chowla, B. Dwork and R. J. Evans
\cite{CDE} in 1986 who used Gauss and Jacobi sums to show that
$$\bi{(p-1)/2}{(p-1)/4}\eq\f{2^{p-1}+1}2\l(2x-\f p{2x}\r)\pmod {p^2}.$$
For more such congruences involving products of one or more binomial
coefficients, the reader may consult the excellent survey \cite{HW}
by R. H. Hudson and K. S. Williams. In 2009 the author (cf. \cite[Conjecture 5.5]{Su2})
conjectured that
$$\sum_{k=0}^{p-1}\f{\bi{2k}k^2}{8^k}\eq\sum_{k=0}^{p-1}\f{\bi{2k}k^2}{(-16)^k}
\eq\l(\f 2p\r)\sum_{k=0}^{p-1}\f{\bi{2k}k^2}{32^k}\eq\l(\f2p\r)\l(2x-\f p{2x}\r)\ (\mo\ p^2),$$
and this was later confirmed by Z.-H. Sun \cite{S1}. Recently the author \cite{Su5} determined
$x$ mod $p^2$ via the congruence
$$\l(\f 2p\r)x\eq\sum_{k=0}^{(p-1)/2}\f{k+1}{8^k}\bi{2k}k^2\eq\sum_{k=0}^{(p-1)/2}\f{2k+1}{(-16)^k}\bi{2k}k^2\ \ (\mo\ p^2).$$
Note that $p\mid\bi{2k}k$ for all $k=(p+1)/2,\ldots,p-1$.

Let $p$ be an odd prime. By \cite{I,vH},
$$\sum_{k=0}^{p-1}\f{\bi{2k}k^3}{64^k}\eq\begin{cases}4x^2-2p\ (\mo\ p^2)&\t{if}\ p=x^2+y^2\ (2\nmid x\ \&\ 2\mid y),
\\0\ (\mo\ p^2)&\t{if}\ p\eq3\ (\mo\ 4).\end{cases}$$
In \cite{Su2,Su4} the author made conjectures on $\sum_{k=0}^{p-1}\bi{2k}k^3/m^k$ mod $p^2$
for $m=1,-8,16,-64,256,-512,4096$; for example,
Conjecture 5.3 of \cite{Su2} states that
$$\sum_{k=0}^{p-1}\bi{2k}k^3\eq\begin{cases}4x^2-2p\ (\mo\ p^2)&\t{if}\ (\f
p7)=1\ \&\ p=x^2+7y^2,\\0\ (\mo\ p^2)&\t{if}\ (\f
p7)=-1.\end{cases}$$
(Throughout this paper, when we write a multiple of a prime in
the form $ax^2+by^2$, we always assume that $x$ and $y$ are {\it nonzero}
integers.)
To attack such conjectures, Z.-H. Sun
\cite{S2} deduced the useful combinatorial identity
\begin{equation}\label{1.1}\sum_{k=0}^n\bi{n+k}{2k}\bi{2k}k^2x^k=P_n(\sqrt{1+4x})^2\end{equation}
where $P_n(x)$ is the Legendre polynomial of degree $n$ given by
$$P_n(x):=\sum_{k=0}^n\bi nk\bi{n+k}k\l(\f{x-1}2\r)^k.$$
Actually (\ref{1.1}) is just a special case of the well-known Clausen formula for hypergeometric series.
We can rewrite (\ref{1.1}) in the form
\begin{equation}\label{1.2}\sum_{k=0}^n\bi nk\bi{n+k}k\bi{2k}k(x(x+1))^k=D_n(x)^2\end{equation}
 where $D_n(x)$ is
the Delannoy polynomial of degree $n$ given by
$$D_n(x):=\sum_{k=0}^n\bi nk\bi{n+k}kx^k.$$
Note that those $D_n=D_n(1)$ $(n=0,1,2,\ldots)$ are central Delannoy
numbers (see, e.g., \cite{CHV}, \cite{Su3} and \cite[p.\,178]{St}.
It is well known that $P_n(-x)=(-1)^nP_n(x)$, i.e.,
$(-1)^nD_n(x)=D_n(-x-1)$ (cf. \cite[Remark 1.2]{Su3}). As observed
by Z.-H. Sun \cite[Lemma 2.2]{S1}, if $0\ls k\ls n=(p-1)/2$ then
$$\bi{n+k}{2k}\eq\f{\bi{2k}k}{(-16)^k}\pmod{p^2}$$
and hence
$$\bi nk\bi{n+k}k=\bi{n+k}{2k}\bi{2k}k\eq\f{\bi{2k}k^2}{(-16)^k}\pmod{p^2}.$$
This simple trick was also realized by van Hamme \cite[p.\,231]{vH}.
Combining this useful trick with the identity (2), we see that
\begin{equation}\label{1.3}\sum_{k=0}^{p-1}\f{\bi{2k}k^3}{(-16)^k}(x(x+1))^k\eq\(\sum_{k=0}^{p-1}\f{\bi{2k}k^2}{(-16)^k}x^k\)^2\pmod{p^2}.\end{equation}
To study the author's conjectures on
$$\sum_{k=0}^{p-1}\f{\bi{2k}k^2\bi{3k}k}{m^k},\ \sum_{k=0}^{p-1}\f{\bi{2k}k^2\bi{4k}{2k}}{m^k},\ \sum_{k=0}^{p-1}\f{\bi{2k}{k}\bi{3k}k\bi{6k}{3k}}{m^k}$$
modulo $p^2$ (with $m$ suitable integers not divisible by $p$) given in \cite{Su4, Su7}, Z.-H. Sun \cite{S2,S3,S4}
managed to prove the following congruences similar to (3):
\begin{align}\sum_{k=0}^{p-1}\f{\bi{2k}k^2\bi{4k}{2k}}{(-64)^k}(x(x+1))^k\eq&\(\sum_{k=0}^{p-1}\f{\bi{4k}{2k}\bi{2k}k}{(-64)^k}x^k\)^2\pmod{p^2},
\\\sum_{k=0}^{p-1}\f{\bi{2k}k^2\bi{3k}k}{(-27)^k}(x(x+1))^k\eq&\(\sum_{k=0}^{p-1}\f{\bi{2k}k\bi{3k}k}{(-27)^k}x^k\)^2\pmod{p^2}\ (p>3),
\\\sum_{k=0}^{p-1}\f{\bi{2k}{k}\bi{3k}k\bi{6k}{3k}}{(-432)^k}(x(x+1))^k\eq&\(\sum_{k=0}^{p-1}\f{\bi{6k}{3k}\bi{3k}k}{(-432)^k}x^k\)^2\pmod{p^2}\ (p>3).
\end{align}

In 1859 G. Bauer proved that
$$\sum_{k=0}^\infty(4k+1)\f{\bi{2k}k^3}{(-64)^k}=\f 2{\pi}.$$
In 1914 S. Ramanujan \cite{R} found 16 new series for $1/\pi$ which are quite similar to Bauer's series.
The (rational) Ramanujan-type series for
$1/\pi$ (cf. B. C. Berndt \cite[pp. 353-354]{Be}, and also \cite{BB} and \cite{ChCh}) have the following form:
\begin{equation}\label{1.7}\sum_{k=0}^\infty(a+dk)\f{f(k)}{m^k}=\f C{\pi},\end{equation}
where $f(k)$ refers to one of
$$\bi{2k}k^3,\ \bi{2k}k^2\bi{3k}k,\ \bi{2k}k^2\bi{4k}{2k},\ \bi{2k}k\bi{3k}k\bi{6k}k,$$
and $a,d,m$ are integers with $dm\not=0$, and $C^2$ is rational.
Up to now, 36 such series have been established via the theory of modular forms.
The reader may consult [CCL], [CC], and S. Cooper [C] for some other series for $1/\pi$.

Let $p$ be an odd prime. Note that $\Gamma(1/2)^2=\pi$ while $\Gamma_p(1/2)^2=(-1)^{(p+1)/2}$,
where $\Gamma_p(x)$ denotes the $p$-adic $\Gamma$-function.
In view of this, in 1997 van Hamme \cite{vH} studied $p$-adic supercongruences
for partial sums of some hypergeometric series involving the Gamma function.
(If a $p$-adic congruence happens to be true modulo a higher power of $p$, then it is called a supercongruence.)
For example, Bauer's series led him to conjecture that
$$\sum_{k=0}^{p-1}(4k+1)\f{\bi{2k}k^3}{(-64)^k}\eq p\l(\f{-1}p\r)\ \ (\mo\ p^3),$$
 which was later confirmed by E. Mortenson \cite{M2} in 2008.
More supercongruences motivated by Ramanujan-type series have been investigated by some followers of van Hamme,
see, e.g., L. Long \cite{L} and Zudilin's work stated there.
The author \cite{Su6} refined the congruence by van Hamme and Mortenson to the following congruence mod $p^4$:
 $$\sum_{k=0}^{p-1}(4k+1)\f{\bi{2k}k^3}{(-64)^k}\eq\sum_{k=0}^{(p-1)/2}(4k+1)\f{\bi{2k}k^3}{(-64)^k}\eq p\l(\f{-1}p\r)+p^3E_{p-3}\ (\mo\ p^4),$$
 where $E_0,E_1,\ldots$ are the Euler numbers defined by
 $$E_0=1,\ \t{and}\ \sum^n_{k=0\atop 2\mid k}\bi nkE_{n-k}=0\ \ (n\in\Z^+).$$
For more conjectural connections between Ramanujan-type congruences
and Euler numbers or Euler polynomials, the reader may consult \cite{Su4}.

Gosper announced in 1974 that $\sum_{k=0}^\infty(25k-3)/(2^k\bi{3k}k)=\pi/2$ (see G. Almkvist, C. Krattenthaler and J. Petersson
\cite{AKP} for a simple proof).
Though this is not a Ramanujan-type series, the author conjectures that for any prime $p>3$ we have
\begin{align*}
2p\sum_{k=0}^{p-1}\f{25k-3}{2^k\bi{3k}k}&\eq3\l(\f{-1}p\r)+(E_{p-3}-9)p^2\ (\mo\ p^3),
\\p\sum_{k=0}^{(p-1)/2}\f{25k-3}{2^k\bi{3k}k}&\eq \l(\f{-1}p\r)-\l(\f 2p\r)\f {5p}2\ (\mo\ p^2),
\\\sum_{k=0}^{p-1}(25k+3)k2^k\bi{3k}k&\eq 6\l(\f{-1}p\r)-18p\ (\mo\ p^2).
\end{align*}
The author \cite{Su4} found some new series for powers of $\pi$ motivated by corresponding $p$-adic congruences.
Here is a new example: Immediately after the author discovered the conjectural congruence
\begin{align*}&\sum_{n=0}^{p-1}\f{28n+5}{576^n}\bi{2n}n\sum_{k=0}^n\f{5^k\bi{2k}k^2\bi{2(n-k)}{n-k}^2}{\bi nk}
\\& \qquad\eq p\l(\f{-1}p\r)\l(3+2\l(\f 2p\r)\r)\pmod{p^2}
\end{align*} for any prime $p>3$, he conjectured (on Jan. 14, 2012) that
\begin{equation}\label{1.8}
\sum_{n=0}^{\infty}\f{28n+5}{576^n}\bi{2n}n\sum_{k=0}^n\f{5^k\bi{2k}k^2\bi{2(n-k)}{n-k}^2}{\bi nk}
=\f 9{\pi}(2+\sqrt2).\end{equation}

The author \cite{Su4,Su7} found that an identity like (\ref{1.7})
usually corresponds to a congruence for $\sum_{k=0}^{p-1}f(k)/m^k$ modulo $p^2$ in terms of parameters in  representations
of a prime $p$ or its multiple by certain binary quadratic forms. This was the main starting point of the author's discoveries of
many new series for $1/\pi$.

Let $n\in\N$. Clearly $\bi{2n}n$ is the coefficient
of $x^n$ in the expansion of $(x^2+2x+1)^n=(x+1)^{2n}$. The $n$th
central trinomial coefficient
$$T_n=[x^n](x^2+x+1)^n$$ is the coefficient of $x^n$ in the
expansion of $(x^2+x+1)^n$. Since $T_n$ is the constant term of
$(1+x+x^{-1})^n$, by the multi-nomial theorem we see that
$$T_n=\sum_{k=0}^{\lfloor n/2\rfloor}\f{n!}{k!k!(n-2k)!}
=\sum_{k=0}^{\lfloor n/2\rfloor}\bi n{2k}\bi{2k}k=\sum_{k=0}^n\bi nk\bi{n-k}k.$$
Central trinomial coefficients arise naturally in enumerative
combinatorics (cf. \cite{Sl}), e.g., $T_n$ is the number of lattice paths
from the point $(0, 0)$ to $(n, 0)$ with only allowed steps $(1,
1)$, $(1, -1)$ and $(1, 0)$.

Given $b,c\in\Z$,  we define the {\it generalized central trinomial
coefficients}
\begin{equation}\begin{aligned}T_n(b,c):=&[x^n](x^2+bx+c)^n=[x^0](b+x+cx^{-1})^n
\\=&\sum_{k=0}^{\lfloor n/2\rfloor}\bi n{2k}\bi{2k}kb^{n-2k}c^k
=\sum_{k=0}^{\lfloor n/2\rfloor}\bi nk\bi{n-k}k b^{n-2k}c^k.
\end{aligned}\label{1.9}\end{equation}
Clearly $T_n(2,1)=\bi{2n}n$ and $T_n(1,1)=T_n$. An efficient way to
compute $T_n(b,c)$ is to use the initial values $T_0(b,c)=1$ and
$T_1(b,c)=b$, and the recursion
$$(n+1)T_{n+1}(b,c)=(2n+1)bT_n(b,c)-n(b^2-4c)T_{n-1}(b,c)\ \ (n=1,2,\ldots).$$
Note that the recursion is rather simple if $b^2-4c=0$.

Let $b,c\in\Z$ and $d=b^2-4c$. It is known that $T_n(b,c)=\sqrt
d^nP_n(b/\sqrt d)$ if $d\not=0$ (see, e.g., \cite{N} and
\cite{Su9}). Thus
\begin{equation}\label{1.10}T_n(b,c)=\sum_{k=0}^n\bi{n+k}{2k}\bi{2k}k\l(\f{b-\sqrt d}2\r)^k\sqrt d^{n-k}.\end{equation}
(In the case $d=0$, (10) holds trivially since
$x^2+bx+c=(x+b/2)^2$.) By the Laplace-Heine formula (cf. \cite[p.\,194]{Sz}), for any complex number $x\not\in[-1,1]$ we have
$$P_n(x)\sim\f{(x+\sqrt{x^2-1})^{n+1/2}}{\sqrt{2n\pi}\root 4\of{x^2-1}}\quad \ \t{as}\ n\to+\infty.$$
It follows that if $b>0$ and $c>0$ then
\begin{equation}\label{1.11}T_n(b,c)\sim f_n(b,c):=\f{(b+2\sqrt{c})^{n+1/2}}{2\root{4}\of c\sqrt{n\pi}}\quad\t{as}\ n\to+\infty.\end{equation}
Note that $T_n(-b,c)=(-1)^nT_n(b,c)$.

The generalized central trinomial coefficients seem to be natural
extensions of the central binomial coefficients.
To see this, in the next section we study congruences for
$$\sum_{k=0}^{p-1}\f{\bi{2k}kT_k(b,c)}{m^k}\ \ \t{and}\ \
\sum_{k=0}^{p-1}\f{\bi{2k}kT_{2k}(b,c)}{m^k}$$
modulo an odd prime $p$, where $b,c,m\in\Z$
and $m\not\eq0\ (\mo\ p)$. One may compare them with congruences for $\sum_{k=0}^{p-1}\bi{2k}k^2/m^k$ mod $p^2$ with $m=8,-16,32$.
Since
$$T_k(2,1)=\bi{2k}k,\ T_{2k}(2,1)=\bi{4k}{2k}\ \t{and}\
T_{3k}(2,1)=\bi{6k}{3k},$$ in Section 3 we are going to investigate general
sums
$$\sum_{k=0}^{p-1}\f{\bi{2k}k^2}{m^k}T_k(b,c),\ \sum_{k=0}^{p-1}\f{\bi{2k}k\bi{3k}k}{m^k}T_k(b,c),
\ \sum_{k=0}^{p-1}\f{\bi{4k}{2k}\bi{2k}k}{m^k}T_k(b,c)$$
and
$$\sum_{k=0}^{p-1}\f{\bi{2k}k^2}{m^k}T_{2k}(b,c),\ \sum_{k=0}^{p-1}\f{\bi{2k}k\bi{3k}k}{m^k}T_{3k}(b,c)$$
modulo $p^2$, where $p$ is an odd prime, $b,c,m\in\Z$ and
$m\not\eq0\ (\mo\ p)$. For this purpose, we need to extend those
congruences (3)-(6) in Section 3.

 Section 4 contains
many conjectural congruences involving generalized central binomial coefficients and they offer backgrounds for those
conjectural series for $1/\pi$ in Sect. 5. In the fifth section we first
show a theorem on dualities and then propose 61 new conjectural series
for $1/\pi$ based on our investigation of congruences.

\section{On $\sum_{k=0}^{p-1}\bi{2k}kT_k(b,c)/m^k$ and $\sum_{k=0}^{p-1}\bi{2k}kT_{2k}(b,c)/m^k$ modulo $p$}
\setcounter{lemma}{0} \setcounter{theorem}{0}
\setcounter{corollary}{0} \setcounter{remark}{0}
\setcounter{conjecture}{0}

\begin{lemma}\label{Lem2.1} Let $p=2n+1$ be an odd prime and let $k\in\{0,\ldots,n\}$. Then
\begin{equation}\label{2.1}\bi{2k}k\eq(-1)^n16^k\bi{2(n-k)}{n-k}\pmod p.\end{equation}
Given $b,c\in\Z$ with $b^2\not\eq 4c\pmod p$, we also have
\begin{equation}\label{2.2} T_{2(n-k)}(b,c)\eq\l(\f{b^2-4c}p\r)\f{T_{2k}(b,c)}{(b^2-4c)^{2k}}\pmod p.\end{equation}
\end{lemma}
\Proof. (12) holds because
$$\f{\bi{2k}k}{(-4)^k}=\bi{-1/2}k\eq\bi nk=\bi n{n-k}\eq \bi{-1/2}{n-k}=\f{\bi{2(n-k)}{n-k}}{(-4)^{n-k}}\pmod p.$$
For $b,c\in\Z$ with $d=b^2-4c\not\eq0\ (\mo\ p)$, we get (13) from
the known result $d^jT_{p-1-j}(b,c)\eq(\f dp)T_j(b,c)$ for
$j=0,\ldots,p-1$ (see \cite[(14)]{N} or \cite[Lemma 2.2]{Su9}). \qed

\begin{theorem}\label{Th2.1} Let $p$ be an odd prime and let $m,b,c\in\Z$ with $m\not\eq0\pmod p$.
If $m\eq4b\pmod p$, then
\begin{equation}\label{2.3}\begin{aligned} &\sum_{k=0}^{p-1}\f{\bi{2k}k}{m^k}T_k(b,c)
\\\eq&\begin{cases}(\f mp)2xc^{(p-1)/4}\pmod p&\t{if}\ p=x^2+y^2\ (4\mid x-1),
\\0\pmod p&\t{if}\ p\eq3\pmod4.\end{cases}
\end{aligned}\end{equation}
If $m\not\eq 4b\pmod{p}$, then
\begin{equation}\label{2.4}\sum_{k=0}^{p-1}\f{\bi{2k}k}{m^k}T_k(b,c)
\eq\l(\f{m(m-4b)}p\r)\sum_{k=0}^{p-1}\f{\bi{4k}{2k}\bi{2k}kc^k}{(m-4b)^{2k}}\pmod p.\end{equation}
Also, provided that $d=b^2-4c\not\eq0\pmod p$, for any $h\in\Z^+$ we have
\begin{equation}\label{2.5}\sum_{k=0}^{p-1}\f{\bi{2k}k^hT_{2k}(b,c)}{m^k}
\eq\l(\f{(-1)^hdm}p\r)\sum_{k=0}^{p-1}\f{\bi{2k}k^hT_{2k}(b,c)}{(16^hd^2/m)^k}\pmod p.\end{equation}
\end{theorem}
\Proof. Set $n=(p-1)/2$. As
$\bi nk\eq\bi{-1/2}k=\bi{2k}k/(-4)^k\pmod p$ for all $k=0,\ldots,p-1$,
we have
\begin{align*} \sum_{k=0}^{p-1}\f{\bi{2k}k}{m^k}T_k(b,c)\eq&\sum_{k=0}^n\bi nk\l(-\f 4m\r)^k[x^0](x+b+cx^{-1})^k
\\=&[x^0]\l(1-\f 4m\cdot\f{x^2+bx+c}x\r)^n
\\\eq&\l(\f mp\r)[x^n](mx-4(x^2+bx+c))^n
\\\eq&(-1)^n\l(\f mp\r)[x^n]\l(x^2-\f{m-4b}4x+c\r)^n
\\=&\l(\f mp\r)T_n\l(\f{m-4b}4,c\r)=\l(\f mp\r)\f{T_n(m-4b,16c)}{2^{2n}}
\\\eq&\l(\f mp\r)T_n(m-4b,16c)\pmod p.
\end{align*}
Observe that
\begin{align*} T_n(m-4b,16c)=&\sum_{k=0}^{\lfloor n/2\rfloor}\bi n{2k}\bi{2k}k(m-4b)^{n-2k}(16c)^k
\\\eq&\sum_{k=0}^{\lfloor n/2\rfloor}\f{\bi{4k}{2k}}{(-4)^{2k}}\bi{2k}k(m-4b)^{n-2k}(16c)^k
\\=&\sum_{k=0}^{\lfloor n/2\rfloor}\bi{4k}{2k}\bi{2k}k(m-4b)^{n-2k}c^k\pmod p.
\end{align*}
Thus (\ref{2.4}) holds when $m\not\eq 4b\pmod p$.
If $m\eq 4b\pmod p$, then
$$ T_n(m-4b,16c)\eq\begin{cases}\bi{2n}n\bi n{n/2}c^{n/2}\pmod p&\t{if}\ 2\mid n,
\\0\pmod p&\t{if}\ 2\nmid n.\end{cases}$$
Clearly $\bi{2n}n=\bi{p-1}n\eq(-1)^n\pmod p$. If $p=2n+1\eq1\pmod4$ and $p=x^2+y^2$ with $x\eq1\pmod4$,
then $\bi{n}{n/2}\eq2x\pmod p$ as observed by Gauss. Thus, (\ref{2.3}) holds when $m\eq 4b\pmod p$.

Now suppose that $d=b^2-4c\not\eq0\pmod p$ and $h\in\Z^+$. In view of Lemma 2.1, we have
\begin{align*} \sum_{k=0}^n\f{\bi{2k}k^hT_{2k}(b,c)}{m^k}
\eq&\sum_{k=0}^n\f{((-1)^n16^k\bi{2(n-k)}{n-k})^h}{m^k}\l(\f dp\r)d^{2k}T_{2(n-k)}(b,c)
\\=&(-1)^{hn}\l(\f dp\r)\sum_{j=0}^n\l(\f{16^hd^2}m\r)^{n-j}\bi{2j}j^hT_{2j}(b,c)
\\\eq&\l(\f{(-1)^hdm}p\r)\sum_{k=0}^n\f{\bi{2k}k^hT_{2k}(b,c)}{(16^hd^2/m)^k}\pmod p.
\end{align*}
Recall that $p\mid\bi{2k}k$ for each $k=n+1,\ldots,p-1$. So (\ref{2.5}) follows.  \qed

\begin{corollary}\label{Cor2.1} Let $p$ be an odd prime. Then
\begin{align*}&\sum_{k=0}^{p-1}\f{\bi{2k}k}{4^k}T_k(1,2)
\\\eq&\begin{cases}(-1)^{(x-1)/2+y/4}2x\ (\mo\ p)&\t{if}\ 8\mid p-1\ \&\ p=x^2+y^2\ (2\nmid x),
\\(-1)^{(y-2)/4}2y\ (\mo\ p)&\t{if}\ 8\mid p-5\ \&\ p=x^2+y^2\ (2\mid y),\\0\pmod{p}&\t{if}\ p\eq3\pmod4.\end{cases}
\end{align*}
\end{corollary}
\Proof. If $p\eq3\pmod4$, then
$\sum_{k=0}^{p-1}\bi{2k}kT_k(1,2)/4^k\eq0\pmod p$ by (\ref{2.3}) with $m=4$, $b=1$ and $c=2$.

Now assume that $p\eq1\pmod4$ and write $p=x^2+y^2$ with $x\eq1\pmod 4$ and $y\eq0\pmod2$.
Applying (\ref{2.3}) with $m=4$, $b=1$ and $c=2$, we get
$$\sum_{k=0}^{p-1}\f{\bi{2k}k}{4^k}T_k(1,2)\eq 2x\times 2^{(p-1)/4}\pmod p.$$
By Exercise 27 of \cite[p.\,64]{IR} (an observation of Dirichlet),
$$2^{(p-1)/4}\eq\l(\f yx\r)^{xy/2}\pmod p.$$
Note that
$$\l(\f yx\r)^2=\f{y^2}{x^2}\eq-1\pmod p\ \t{and hence}\ \l(\f yx\r)^4\eq1\pmod 4.$$
So we have
$$2^{(p-1)/4}\eq\l(\f yx\r)^{y/2}\pmod p.$$
If $p\eq1\pmod 8$, then $4\mid y$ and hence $2^{(p-1)/4}\eq(-1)^{y/4}\pmod p$.
If $p\eq5\pmod 8$, then $y\eq2\pmod 4$ and hence
$$2^{(p-1)/4}\eq\l(\f yx\r)^{2(y-2)/4}\f yx\eq(-1)^{(y-2)/4}\ \f yx\pmod p.$$

Combining the above, we obtain the desired result. \qed

\begin{corollary}\label{Cor2.2} For any prime $p>3$ we have
$$\sum_{k=0}^{p-1}\f{\bi{2k}kT_k}{(-4)^k}\eq\l(\f{-1}p\r)\pmod p
\ \ \t{and}\ \ \sum_{k=0}^{p-1}\f{\bi{2k}kT_k}{12^k}\eq\l(\f p3\r)\pmod p.$$
\end{corollary}
\Proof. Applying (\ref{2.4}) with $b=c=1$ and $m\in\{-4,12\}$ we obtain
$$\sum_{k=0}^{p-1}\f{\bi{2k}kT_k}{(-4)^k}\eq\l(\f{(-4)(-8)}p\r)\sum_{k=0}^{p-1}\f{\bi{4k}{2k}\bi{2k}k}{64^k}\pmod p$$
and
$$\sum_{k=0}^{p-1}\f{\bi{2k}kT_k}{12^k}\eq\l(\f{12\times 8}p\r)\sum_{k=0}^{p-1}\f{\bi{4k}{2k}\bi{2k}k}{64^k}\pmod p.$$
It is known that
$$\sum_{k=0}^{p-1}\f{\bi{4k}{2k}\bi{2k}k}{64^k}\eq\l(\f{-2}p\r)\pmod {p^2},$$
which was conjectured in \cite{RV} and proved in \cite{M1}.
So the two congruences in Corollary 2.2 are valid. \qed

\begin{theorem}\label{Thm2.2} Let $p$ be an odd prime. Then
$$\sum_{k=0}^{p-1}\f{\bi{2k}kT_{2k}}{4^k}\eq\l(\f{-2}p\r)\pmod p.$$
\end{theorem}
\Proof. Set $n=(p-1)/2$. Then
\begin{align*} \sum_{k=0}^{p-1}\f{\bi{2k}k}{4^k}T_{2k}
\eq&\sum_{k=0}^n\bi nk(-1)^k[x^0](1+x+x^{-1})^{2k}
\\=&[x^0]\l(1-(1+x+x^{-1})^2\r)^n
\\=&[x^0](-1)^n\l(\f{x^2+1}x\cdot\f{(x+1)^2}x\r)^n
\\=&(-1)^n[x^{2n}](x^2+1)^n(x+1)^{2n}=(-1)^n\sum_{k=0}^n\bi nk\bi{2n}{2k}
\\\eq&\sum_{k=0}^n\bi nk(-1)^n=(-2)^n\eq\l(\f{-2}p\r)\pmod p.
\end{align*}
This concludes the proof. \qed

\begin{remark}\label{Rem2.1} For any prime $p>3$ we observe the following congruences:
\begin{align*}\sum_{k=0}^{p-1}\f{\bi{2k}k}{4^k}T_{2k}(5,4)\eq&1\pmod{p},
\ \ \sum_{k=0}^{p-1}\f{\bi{2k}k}{4^k}T_{2k}(3,1)\eq\l(\f2p\r)\pmod p,
\\\sum_{k=0}^{p-1}\f{\bi{2k}k}{16^k}T_{2k}(4,9)\eq&\l(\f p3\r)\,(\mo\ p),\,
 \sum_{k=0}^{p-1}\f{\bi{2k}k}{16^k}T_{2k}(8,25)\eq\l(\f{-5}p\r)\, (\mo\ p).
\end{align*}
\end{remark}

\begin{conjecture}\label{Conj2.1} Let $p>3$ be a prime. Then
\begin{align*}\sum_{k=0}^{p-1}\f{\bi{2k}k}{12^k}T_k\eq&\l(\f p3\r)\f{3^{p-1}+3}4\pmod{p^2},
\\\sum_{k=0}^{(p-1)/2}\f{\bi{2k}k}{16^k}T_{2k}(4,1)\eq&1\pmod{p^2},
\\\sum_{k=0}^{(p-1)/2}\f{\bi{2k}k}{(k+1)16^k}T_{2k}(4,1)\eq&\f 43\l(\l(\f 3p\r)-p\l(\f{-1}p\r)\r)\pmod{p^2},
\\\sum_{k=0}^{(p-1)/2}\f{\bi{2k}k}{4^k}T_{2k}(3,4)\eq&\l(\f{-1}p\r)\f{7-3^p}4\pmod{p^2},
\\\sum_{k=0}^{(p-1)/2}\f{\bi{2k}k}{16^k}T_{2k}(8,9)\eq&\l(\f 3p\r)\pmod{p^2},
\end{align*}
and
$$\sum_{k=0}^{p-1}\f{\bi{3k}k}{432^k}T_{3k}(6,1)\eq1\pmod{p}.$$
\end{conjecture}

\begin{conjecture}\label{Conj2.2} Let $p>3$ be a prime. Then
\begin{align*}&\sum_{k=0}^{(p-1)/2}\f{\bi{2k}k}{16^k}T_{2k}(2,3)\eq\sum_{k=0}^{(p-1)/2}\f{\bi{2k}k}{16^k}T_{2k}(4,-3)
\\\eq&\begin{cases}(\f{-1}p)(\f x3)(2x-\f p{2x})\pmod{p^2}&\t{if}\ p=x^2+3y^2,
\\0\pmod p&\t{if}\ p\eq2\pmod3.\end{cases}
\end{align*}
Also,
$$\sum_{k=0}^{(p-1)/2}\f{\bi{2k}k}{4^k}T_{2k}(1,-3)\eq\begin{cases}(-1)^{xy/2}(\f x3)2x\pmod p&\t{if}\ p=x^2+3y^2,
\\0\pmod{p}&\t{if}\ p\eq2\,(\mo\ 3);\end{cases}$$
and
\begin{align*}&\sum_{k=0}^{(p-1)/2}\f{\bi{2k}k}{16^k}T_{2k}(4,3)
\\\eq&\begin{cases}(-1)^{\lfloor x/6\rfloor+y/2}2x\,(\mo\ p)&\t{if}\ 12\mid p-1\, \&\,p=x^2+y^2\, (2\nmid x),
\\(-1)^{(x+y+1)/2}(\f{xy}3)2y\,(\mo\ p)&\t{if}\ 12\mid p-5\,\&\, p=x^2+y^2\, (2\nmid x),
\\0\,(\mo\ p)&\t{if}\ p\eq3\pmod4.\end{cases}
\end{align*}
\end{conjecture}

\begin{conjecture}\label{Conj2.3} Let $p$ be an odd prime. Then
\begin{align*}&\sum_{k=0}^{(p-1)/2}\f{\bi{2k}k}{16^k}T_{2k}(12,-7)
\\\eq&\begin{cases}2x(\f x7)\pmod p&\t{if}\ (\f p7)=1\ \&\ p=x^2+7y^2,
\\0\pmod{p}&\t{if}\ (\f p7)=-1,\ \t{i.e.},\ p\eq3,5,6\pmod7.\end{cases}
\end{align*}
\end{conjecture}

\section{Extensions of (2)--(6) with applications to sums involving generalized central trinomial coefficients}
\setcounter{lemma}{0} \setcounter{theorem}{0}
\setcounter{corollary}{0} \setcounter{remark}{0}
\setcounter{conjecture}{0}

 Our following theorem is a natural generalization of (2).

\begin{theorem}\label{Th3.1} For any $n\in\N$ we have
\begin{equation}\label{3.1}D_n(x)D_n(y)=\sum_{k=0}^n\bi{n+k}{2k}\bi{2k}k\sum_{j=0}^k\bi{k+j}{2j}\bi{2j}j(xy+y)^j(x-y)^{k-j},\end{equation}
\end{theorem}
\Proof.  Let $a_n$ denote the left hand side
or the right-hand side of (17). It is easy to see that
$$a_0=1,\ a_1=(2x+1)(2y+1),\ a_2=(6x^2+6x+1)(6y^2+6y+1)$$
and
$$a_3=(20x^3+30x^2+12x+1)(20y^3+30y^2+12y+1).$$
Applying the Zeilberger algorithm (cf. \cite[pp.\,101-119]{PWZ}) via {\tt Mathematica} we find the recursion for $n\gs 3$:
\begin{align*}&(n+1)^2(2n-3)a_{n+1}-(2n-3)(2n+1)^2(2x+1)(2y+1)a_n
\\&+(2n-1)A(n,x,y)a_{n-1}-(2n-3)^2(2n+1)(2x+1)(2y+1)a_{n-2}
\\&+(n-2)^2(2n+1)a_{n-3}
\\&=0,
\end{align*}
where
$$A(n,x,y):=6n^2-6n-5+(16n^2-16n-12)(x+y-x^2-y^2).$$
Thus (17) holds by induction. \qed

Now we give our extensions of $(3)-(6)$.

\begin{theorem}\label{Th3.2} Let $p$ be a prime and let
$a\in\Z^+$. Let $h$ be a $p$-adic integer and set
$w_k(h)=\bi hk\bi{h+k}k$ for $k\in\N$. Then
\begin{equation}\begin{aligned} &\(\sum_{k=0}^{p^a-1}w_k(h)x^k\)\(\sum_{k=0}^{p-1}w_k(h)y^k\)
\\\eq&\sum_{k=0}^{p^a-1}w_k(h)\sum_{j=0}^k\bi{k+j}{2j}\bi{2j}j(xy+y)^j(x-y)^{k-j}\pmod{p^2}.
\end{aligned}\label{3.2}\end{equation}
In particular, if $p\not=2$ then
\begin{equation}\begin{aligned} &\(\sum_{k=0}^{p^a-1}\f{\bi{2k}k^2}{(-16)^k}x^k\)\(\sum_{k=0}^{p^a-1}\f{\bi{2k}k^2}{(-16)^k}y^k\)
\\\eq&\sum_{k=0}^{p^a-1}\f{\bi{2k}k^2}{(-16)^k}\sum_{j=0}^k\bi{k+j}{2j}\bi{2j}j(xy+y)^j(x-y)^{k-j}\pmod{p^2}
\end{aligned}\label{3.3}\end{equation}
and
\begin{equation}\begin{aligned} &\(\sum_{k=0}^{p^a-1}\f{\bi{4k}{2k}\bi{2k}k}{(-64)^k}x^k\)\(\sum_{k=0}^{p^a-1}\f{\bi{4k}{2k}\bi{2k}k}{(-64)^k}y^k\)
\\\eq&\sum_{k=0}^{p^a-1}\f{\bi{4k}{2k}\bi{2k}k}{(-64)^k}\sum_{j=0}^k\bi{k+j}{2j}\bi{2j}j(xy+y)^j(x-y)^{k-j}\pmod{p^2};
\end{aligned}\label{3.4}\end{equation}
provided $p>3$ we have
\begin{equation}\begin{aligned} &\(\sum_{k=0}^{p^a-1}\f{\bi{2k}k\bi{3k}k}{(-27)^k}x^k\)\(\sum_{k=0}^{p^a-1}\f{\bi{2k}k\bi{3k}k}{(-27)^k}y^k\)
\\\eq&\sum_{k=0}^{p^a-1}\f{\bi{2k}k\bi{3k}k}{(-27)^k}\sum_{j=0}^k\bi{k+j}{2j}\bi{2j}j(xy+y)^j(x-y)^{k-j}\pmod{p^2}
\end{aligned}\label{3.5}\end{equation}
and
\begin{equation}\begin{aligned}&\(\sum_{k=0}^{p^a-1}\f{\bi{6k}{3k}\bi{3k}k}{(-432)^k}x^k\)\(\sum_{k=0}^{p^a-1}\f{\bi{6k}{3k}\bi{3k}k}{(-432)^k}y^k\)
\\\eq&\sum_{k=0}^{p^a-1}\f{\bi{6k}{3k}\bi{3k}k}{(-432)^k}\sum_{j=0}^k\bi{k+j}{2j}\bi{2j}j(xy+y)^j(x-y)^{k-j}\pmod{p^2}.
\end{aligned}\label{3.6}\end{equation}
\end{theorem}
\begin{remark}\label{Rem3.1} Note that
\begin{align*} w_k\l(-\f12\r)=&\f{\bi{2k}k^2}{(-16)^k},\ \ w_k\l(-\f14\r)=\f{\bi{4k}{2k}\bi{2k}k}{(-64)^k},
\\w_k\l(-\f13\r)=&\f{\bi{2k}{k}\bi{3k}k}{(-27)^k},\ \
w_k\l(-\f16\r)=\f{\bi{6k}{3k}\bi{3k}k}{(-432)^k}.\end{align*}
Also, (19)-(22) in the case $x=y$ and $a=1$ yield (3)-(6)
respectively.
\end{remark}

The reader may wonder how we found Theorem 3.2. In fact, (17) is our main clue to the
congruence (19). By refining our proof of (19)-(22) we found (18).

To prove Theorem 3.2 we need two lemmas.

\begin{lemma}\label{Lem3.1} For $m,n\in\N$ we have
\begin{equation}\label{3.7}\sum_{k=0}^n\bi nk\bi{k+m}nw_{k+m}(h)=\f{w_m(h)w_n(h)}{\bi{m+n}n},\end{equation}
where $w_k(h)=\bi hk\bi{h+k}k$ as defined in Theorem $3.2$.
\end{lemma}
\Proof. Let $u_n$ denote the left-hand side of (\ref{3.7}). By applying
the Zeilberger algorithm via {\tt Mathematica}, we find
the recursion:
$$(n+1)(m+n+1)u_{n+1}=(h-n)(h+n+1)u_n\ (n=0,1,2,\ldots).$$
Thus (\ref{3.7}) can be easily proved by induction on $n$. \qed

\begin{lemma}\label{Lem3.2} For $k,m,n\in\N$ we have the combinatorial
identity
\begin{equation}\begin{aligned}&\sum_{j=0}^m(-1)^{m-j}\bi{m+j}{2j}\bi{2j}j\bi{j+k+m}k\bi jn
\\&\quad=\bi{k+m+n}m\bi{k+m}m\bi mn.\end{aligned}\label{3.8}\end{equation}
\end{lemma}
\Proof. If $m<n$ then both sides of (\ref{3.8}) vanish. (\ref{3.8}) in the case
$m=n$ can be directly verified. Let $s_m$ denote the left-hand side
of (\ref{3.8}). By the Zeilberger algorithm we find the recursion
$$(m+1)(m-n+1)s_{m+1}=(k+m+1)(k+m+n+1)s_m\ \ (m=n,n+1,\ldots).$$
So we can show (\ref{3.8}) by induction. \qed

\medskip
\noindent{\it Proof of Theorem} 3.2. In view of Remark 3.1, it
suffices to prove (18). Note that both sides of (18) are
polynomials in $x$ and $y$ and the degrees with respect to $x$ or
$y$ are all smaller than $p^a$.

 Fix $m,n\in\{0,\ldots,p^a-1\}$ and let $c(m,n)$ denote the coefficient of $x^ny^m$ in the right-hand side of (18).
 Define $\bi z{-k}=0$ for $k=1,2,3,\ldots$.
Then $c(m,n)$ coincides with
\begin{align*} &[x^n]\sum_{0\ls j\ls k<p^a}w_k(h)\bi{k+j}{2j}\bi{2j}j(x+1)^j\bi{k-j}{m-j}(-1)^{m-j}x^{k-m}
\\=&\sum_{k=m}^{p^a-1}w_k(h)\sum_{j=0}^m(-1)^{m-j}\bi{k+j}{2j}\bi {2j}j\bi{k-j}{m-j}\bi j{m+n-k}
\\=&\sum_{k=0}^{p^a-1-m}w_{k+m}(h)\sum_{j=0}^m(-1)^{m-j}\bi{k+m+j}{2j}\bi{2j}j\bi{k+m-j}{k}\bi j{n-k}
\\=&\sum_{k=0}^{p^a-1-m}w_{k+m}(h)\sum_{j=0}^m(-1)^{m-j}\bi{m+j}{2j}\bi{2j}j\bi{k+m+j}{k}\bi j{n-k}.
\end{align*}
Applying Lemma 3.2 we get
\begin{align*} c(m,n)=&\bi{m+n}m\sum_{k=0}^{p^a-1-m}w_{k+m}(h)\bi{k+m}m\bi m{n-k}
\\=&\bi{m+n}m\sum_{k=0}^{p^a-1-m}w_{k+m}(h)\bi{k+m}n\bi nk.
\end{align*}
By Lemma 3.1,
\begin{align*}&\sum_{k=0}^{p^a-1}w_{k+m}(h)\bi{k+m}n\bi nk
\\=&\sum_{k=0}^nw_{k+m}(h)\bi{k+m}n\bi nk=\f{w_m(h)w_n(h)}{\bi{m+n}m}.\end{align*}
So, it remains to show
\begin{equation}\label{3.9}\bi{m+n}m\sum_{k=p^a-m}^{p^a-1}w_{k+m}(h)\bi{k+m}n\bi nk\eq0\pmod{p^2}.\end{equation}

To prove (\ref{3.9}) we only need to show
$$\bi{m+n}m\eq\bi{k+m}n\eq0\pmod{p}$$
under the supposition $n\gs k\gs p^a-m$. Note that $m+n\gs k+m\gs
p^a$ and $0<p^a-n\ls k+m-n\ls m<p^a$. As the addition of $m$ and $n$
in base $p$ has at least one carry, we have $p\mid \bi{m+n}m$ by
Kummer's theorem (cf. \cite[p.\,24]{Ri}). Similarly, $p\mid\bi{k+m}n$.

So far we have completed the proof of Theorem 3.2. \qed

Theorem 3.2 implies the following useful result on congruences for
sums of central binomial coefficients and generalized central
trinomial coefficients.

\begin{theorem}\label{Th3.3} Let $p$ be an odd prime and let $x$ be a
$p$-adic integer.
 Let $a\in\Z^+$, $b,c\in\Z$ and $d=b^2-4c$. Set $D:=1+2bx+dx^2$. Then we have
\begin{equation}\begin{aligned}&\sum_{k=0}^{p^a-1}\f{\bi {2k}k^2}{(-16)^k}T_k(b,c)x^k
\\\eq&\(\sum_{k=0}^{p^a-1}\f{\bi{2k}k^2}{32^k}(1-\sqrt{D}+\sqrt d\,x)^k\)
\(\sum_{k=0}^{p^a-1}\f{\bi{2k}k^2}{32^k}(1-\sqrt{D}-\sqrt d\,x)^k\)
\\\eq &P_{(p^a-1)/2}(\sqrt D+\sqrt d\,x)P_{(p^a-1)/2}(\sqrt D-\sqrt d\,x)\pmod{p^2}
\end{aligned}\label{3.10}\end{equation}
and
\begin{equation}\begin{aligned}\sum_{k=0}^{p^a-1}\f{\bi{4k}{2k}\bi{2k}k}{(-64)^k}T_k(b,c)x^k
\eq&\(\sum_{k=0}^{p^a-1}\f{\bi{4k}{2k}\bi{2k}k}{128^k}(1-\sqrt{D}+\sqrt d\,x)^k\)
\\&\times\sum_{k=0}^{p^a-1}\f{\bi{4k}{2k}\bi{2k}k}{128^k}(1-\sqrt{D}-\sqrt d\,x)^k\pmod{p^2}.
\end{aligned}\label{3.11}\end{equation}
If $p>3$, then
\begin{equation}\begin{aligned}\sum_{k=0}^{p^a-1}\f{\bi {2k}k\bi{3k}k}{(-27)^k}T_k(b,c)x^k
\eq&\(\sum_{k=0}^{p^a-1}\f{\bi{2k}k\bi{3k}k}{54^k}(1-\sqrt{D}+\sqrt d\,x)^k\)
\\&\times\sum_{k=0}^{p^a-1}\f{\bi{2k}k\bi{3k}k}{54^k}(1-\sqrt{D}-\sqrt d\,x)^k\pmod{p^2}
\end{aligned}\label{3.12}\end{equation}
and
\begin{equation}\begin{aligned}\sum_{k=0}^{p^a-1}\f{\bi {6k}{3k}\bi{3k}k}{(-432)^k}T_k(b,c)x^k
\eq&\(\sum_{k=0}^{p^a-1}\f{\bi{6k}{3k}\bi{3k}k}{864^k}(1-\sqrt{D}+\sqrt d\,x)^k\)
\\&\times\sum_{k=0}^{p^a-1}\f{\bi{6k}k\bi{3k}k}{864^k}(1-\sqrt{D}-\sqrt d\,x)^k\pmod{p^2}.
\end{aligned}\label{3.13}\end{equation}
\end{theorem}
\begin{remark}\label{Rem3.2} Note that $\sqrt d$ and  $\sqrt D$ in Theorem 3.3 are
viewed as algebraic $p$-adic integers.
\end{remark}

\medskip
\noindent{\it Proof of Theorem} 3.3. Let $n=(p^a-1)/2$. For
$k=0,\ldots,n$ we have
$$\bi{n+k}{2k}=\f{\bi{2k}k}{(-16)^k}\prod_{0<j\ls k}\l(1-\f{p^{2a}}{(2j-1)^2}\r)\eq\f{\bi{2k}k}{(-16)^k}\pmod{p^2}$$
and hence
$$\bi nk\bi{n+k}k=\bi{n+k}{2k}\bi{2k}k\eq\f{\bi{2k}k^2}{(-16)^k}\pmod{p^2}.$$
Note also that $p\mid\bi{2k}k$ for $k=n+1,\ldots,p^a-1$ by Kummer's
theorem. Thus
\begin{align*} P_n(t)=&\sum_{k=0}^n\bi nk\bi{n+k}k\l(\f{t-1}2\r)^k
\\\eq&\sum_{k=0}^n\f{\bi{2k}k^2}{(-16)^k}\l(\f{t-1}2\r)^k
\eq\sum_{k=0}^{p^a-1}\f{\bi{2k}k^2}{32^k}(1-t)^k\pmod{p^2},
\end{align*}
and hence the second congruence in (26) follows.

Set
$$u=\f{\sqrt D+\sqrt d\,x-1}2\ \ \t{and}\ \ v=\f{\sqrt D-\sqrt d\,x-1}2.$$
Then
$$uv+v=\f{D-(\sqrt d\,x+1)^2}4=\f{b-\sqrt d}2x\ \ \ \t{and}\ \ \ u-v=\sqrt d\,x.$$
In view of (\ref{1.10}), for any $k\in\N$ we have
\begin{align*}&\sum_{j=0}^k\bi{k+j}{2j}\bi{2j}j(uv+v)^j(u-v)^{k-j}
\\=&x^k\sum_{j=0}^k\bi{k+j}{2j}\bi{2j}j\l(\f{b-\sqrt d}2\r)^j\sqrt d^{k-j}=x^kT_k(b,c).
\end{align*}
So the first congruence in (26) follows from (19). Similarly,
(27)-(29) are consequences of (20)-(22) respectively.  \qed

For $d\in\{2,3,4,7\}$, it is well known that an odd prime $p$ can be
written in the form $x^2+dy^2$ with $x,y\in\Z$ if and only if
$(\f{-d}p)=1$ (see, e.g., \cite{BEW} and \cite{Co}).

Applying (\ref{3.10}) we get the following new results.

\begin{theorem}\label{Th3.4} Let $p$ be an odd prime. Then
\begin{equation}\sum_{k=0}^{p-1}\f{\bi{2k}k^2T_k(1,-2)}{32^k}
\eq\begin{cases}(\f 2p)(4x^2-2p)\ (\mo\ p^2)&\t{if}\  p=x^2+4y^2,
\\0\ (\mo\ p^2)&\t{if}\ p\eq3\ (\mo\ 4).\end{cases}
\label{3.14}\end{equation}
Also,
\begin{equation}\sum_{k=0}^{p-1}\f{\bi{2k}k^2T_k(2,-1)}{8^k}
\eq\begin{cases}(\f{-1}p)4x^2\ (\mo\ p)&\t{if}\ p=x^2+2y^2,
\\0\ (\mo\ p^2)&\t{if}\ p\eq5,7\ (\mo\ 8);\end{cases}
\label{3.15}\end{equation}
\begin{equation}\sum_{k=0}^{p-1}\f{\bi{2k}k^2T_k(4,1)}{(-4)^k}
\eq\begin{cases} 4x^2\pmod{p}&\t{if}\ p=x^2+3y^2,
\\0\pmod{p^2}&\t{if}\ p\eq2\ (\mo\ 3);\end{cases}
\label{3.16}\end{equation}
and
\begin{equation}\begin{aligned}&\sum_{k=0}^{p-1}(-1)^k\bi{2k}k^2T_k(16,1)\eq\l(\f{-1}p\r)\sum_{k=0}^{p-1}\f{\bi{2k}k^2T_k(1,16)}{(-256)^k}
\\\eq&\begin{cases}4x^2\ (\mo\ p)&\t{if}\ (\f p7)=1\ \t{and}\ p=x^2+7y^2,
\\0\ (\mo\ p^2)&\t{if}\ (\f p7)=-1,\ \t{i.e.},\ p\eq3,5,6\ (\mo\ 7).\end{cases}
\end{aligned}\label{3.17}\end{equation}
\end{theorem}

\begin{remark}\label{Rem3.3} Let $p$ be an odd prime. We guess that $4x^2\ (\mo\ p)$ in (31)-(33) can be replaced by $4x^2-2p\ (\mo\ p^2)$.
\end{remark}

To prove Theorem 3.4 we need a lemma.

\begin{lemma}\label{Lem3.3} Let $p$ be an odd prime. Then
\begin{equation}\label{3.18}\sum_{k=0}^{p-1}\f{\bi{2k}k^2}{32^k}x^k\eq \l(\f 2p\r)x^nP_n\l(1-\f 4x\r)\pmod{p}\end{equation}
and
\begin{equation}\label{3.19}P_n(x)\eq(2x+2)^nP_n\l(\f{3-x}{1+x}\r)\pmod{p}.\end{equation}
\end{lemma}
\Proof. With the help of Lemma 2.1, we get
\begin{align*} \sum_{k=0}^{p-1}\f{\bi{2k}k^2}{32^k}x^k\eq&\sum_{k=0}^n\f{\bi{2k}k^2}{32^k}x^k
\eq\sum_{k=0}^n\f{256^k\bi{2(n-k)}{n-k}^2}{32^k}x^k=\sum_{k=0}^n\bi{2k}k^2(8x)^{n-k}
\\\eq&\l(\f8p\r)x^n\sum_{k=0}^n\f{\bi{2k}k^2}{(-16)^k}\l(-\f 2x\r)^k
\\\eq&\l(\f 2p\r)x^n\sum_{k=0}^n\bi nk\bi{n+k}k\l(\f{(1-4/x)-1}2\r)^k
\\=&\l(\f2p\r)x^nP_n\l(1-\f 4x\r)\pmod{p}.
\end{align*}
This proves (\ref{3.18}).

(\ref{3.19}) follows from \cite[Theorem 2.6]{S1} and its proof. \qed
\medskip

\noindent {\it Proof of Theorem} 3.4. For convenience we set
$n=(p-1)/2$.

(i) Applying (\ref{3.10}) with $b=1$, $c=-2$ and $x=-1/2$, we obtain that
$$\sum_{k=0}^{p-1}\f{\bi{2k}k^2}{32^k}T_k(1,-2)\eq\(\sum_{k=0}^{p-1}\f{\bi{2k}k^2}{32^k}\) \(\sum_{k=0}^{p-1}\f{\bi{2k}k^2}{(-16)^k}\)
\pmod{p^2}.$$
The author \cite[Conjecture 5.5]{Su2} conjectured that
$\sum_{k=0}^{p-1}\bi{2k}k^2/32^k\eq0\ (\mo\ p^2)$ if $p\eq3\ (\mo\
4)$, and
$$\sum_{k=0}^{p-1}\f{\bi{2k}k^2}{32^k}\eq\l(\f 2p\r)\sum_{k=0}^{p-1}\f{\bi{2k}k^2}{(-16)^k}\eq2x-\f p{2x}\ (\mo\ p^2)$$
if $p=x^2+y^2$ with $x\eq1\ (\mo\ 4)$ and
$y\eq0\ (\mo\ 2)$. This was confirmed by Z.-H. Sun \cite{S1}. So the
desired (\ref{3.14}) follows.

(ii) Applying (\ref{3.10}) with $b=2$, $c=-1$ and $x=-2$ we get
\begin{align*} \sum_{k=0}^{p-1}\f{\bi{2k}k^2T_k(2,-1)}{8^k}\eq
\sum_{k=0}^n\f{\bi{2k}k^2}{32^k}\al^k\times \sum_{k=0}^n\f{\bi{2k}k^2}{32^k}\beta^k\ (\mo\ p^2).\end{align*}
where $\al=-4(1+\sqrt2)$ and $\beta=-4(1-\sqrt2)$. Clearly
$\al\beta=-16$. By Lemma 3.3,
$$\sum_{k=0}^n\f{\bi{2k}k^2}{32^k}\al^k\eq\l(\f 2p\r)\al^nP_n(\sqrt2)\pmod{p}$$
and
$$\sum_{k=0}^n\f{\bi{2k}k^2}{32^k}\beta^k\eq\l(\f 2p\r)\beta^n P_n(-\sqrt2)=\l(\f{-2}p\r)\beta^nP_n(\sqrt2)\pmod{p}.$$
By \cite[Theorem 2.7]{S2}, $P_n(\sqrt{2})\eq0\ (\mo\ p)$ if $(\f {-2}p)=-1$,
and $P_n(\sqrt{2})^2\eq(\f{-1}p)4x^2\ (\mo\ p)$ if $(\f {-2}p)=1$ and
$p=x^2+2y^2\ (x,y\in\Z)$. So (\ref{3.15}) holds.

(iii) (\ref{3.10}) with $b=x=4$ and $c=1$ yields that
$$\sum_{k=0}^{p-1}\f{\bi{2k}k^2}{(-4)^k}T_k(4,1)\eq P_n(15+8\sqrt3)P_n(15-8\sqrt3)\pmod{p^2}.$$
 By Lemma 3.3,
 \begin{align*}&(\pm1)^nP_n\l(\f{\sqrt3}2\r)=P_n\l(\pm\f{\sqrt3}2\r)
 \\\eq&(2\pm\sqrt3)^nP_n\l(\f{3\mp\sqrt3/2}{1\pm\sqrt3/2}\r)
 =(2\pm\sqrt3)^nP_n(15\mp8\sqrt3)\pmod{p}.
 \end{align*}
 By \cite[Theorem 2.8]{S2}, $P_n(\sqrt3/2)\eq0\ (\mo\ p)$ if $p\eq2\ (\mo\ 3)$, and
 $P_n(\sqrt3/2)^2\eq(-1)^n4x^2\ (\mo\ p)$ if $p\eq1\ (\mo\ 3)$ and $p=x^2+3y^2\ (x,y\in\Z)$.
 Therefore (\ref{3.16}) is valid.

(iv) Applying (\ref{3.10}) with $b=1$, $c=16$ and $x=1/16$ we obtain that
$$\sum_{k=0}^{p-1}\f{\bi{2k}k^2T_k(1,16)}{(-256)^k}
\eq\(\sum_{k=0}^n\f{\bi{2k}k^2}{32^k}\al^k\)\times \sum_{k=0}^n\f{\bi{2k}k^2}{32^k}\beta^k
\ (\mo\ p^2),$$
where $\al=(1+3\sqrt{-7})/16$ and $\beta=(1-3\sqrt{-7})/16$. Note
that $\al\beta=1/4$. By Lemma 3.3,
$$\sum_{k=0}^{p-1}\f{\bi{2k}k^2}{32^k}\al^k\eq\l(\f2p\r)\al^nP_n(\sqrt{-63})
\ (\mo\ p)$$
and
$$\sum_{k=0}^{p-1}\f{\bi{2k}k^2}{32^k}\beta^k\eq\l(\f2p\r)\beta^nP_n(-\sqrt{-63})=\l(\f{-2}p\r)\beta^nP_n(\sqrt{-63})\pmod{p}.$$

(\ref{3.10}) with $b=16$, $c=1$ and $x=16$ yields that
$$\sum_{k=0}^{p-1}(-1)^k\bi{2k}k^2T_k(16,1)\eq P_n(255+96\sqrt7)P_n(255-96\sqrt7)\pmod{p^2}.$$
By Lemma 3.3,
$$(\pm1)^nP_n\l(\f{3\sqrt7}8\r)\eq (8\pm 3\sqrt7)^nP_n(255\mp 96\sqrt7)\pmod{p}.$$
Therefore
$$(8^2-9\times7)^n\sum_{k=0}^{p-1}(-1)^k\bi{2k}k^2T_k(16,1)\eq (-1)^nP_n\l(\f{3\sqrt7}8\r)^2\pmod{p}.$$

By \cite[Theorem 2.5]{S2}, $P_n(\sqrt{-63})\eq P_n(3\sqrt7/8)\eq0\ (\mo\
p)$ if $(\f p7)=-1$, and
$$P_n(\sqrt{-63})^2\eq (-1)^n P_n\l(\f {3\sqrt7}8\r)^2\eq4x^2\ (\mo\ p)$$
if $(\f p7)=1$ and $p=x^2+7y^2\ (x,y\in\Z)$. Therefore (\ref{3.17}) holds.
\qed

\medskip

Motivated by Theorem 3.4 and the congruence
$$\sum_{k=0}^{p-1}(21k+8)\bi{2k}k^3\eq 8p+16p^4B_{p-3}\ (\mo\ p^5)$$
proved in \cite{Su4} (where $B_0,B_1,\ldots$ are Bernoulli numbers),
 we conjecture that
\begin{align*}\sum_{k=0}^{p-1}(3k+1)\f{\bi{2k}k^2T_k(1,-2)}{32^k}\eq&\l(\f{-2}p\r)\f{2p}{3-(\f{-1}p)}\pmod{p^2},
\\\sum_{k=0}^{p-1}(5k+2)\f{\bi{2k}k^2T_k(2,-1)}{8^k}\eq& p+p\l(\f{-1}p\r)\pmod{p^2},
\\\sum_{k=0}^{p-1}(5k+2)\f{\bi{2k}k^2T_k(4,1)}{(-4)^k}\eq&\f 23p\l(2\l(\f{-1}p\r)+1\r)\pmod{p^2},
\\\sum_{k=0}^{p-1}(255k+112)(-1)^k\bi{2k}k^2T_k(16,1)\eq& 16p\l(3+4\l(\f{-1}p\r)\r)\pmod{p^2},
\\\sum_{k=0}^{p-1}(30k+7)\f{\bi{2k}k^2T_k(1,16)}{(-256)^k}\eq& 7p\l(\f{-1}p\r)\pmod{p^2}.
\end{align*}
The last congruence led the author to find the conjectural identity
$$\sum_{k=0}^{\infty}\f{30k+7}{(-256)^k}\bi{2k}k^2T_k(1,16)=\f{24}{\pi}$$
in Jan. 2011 which was the starting point of the discovery of many
series for $1/\pi$ of new types given in Section 5.

\section{Conjectural congruences related to representations of primes by binary quadratic forms}
\setcounter{lemma}{0} \setcounter{theorem}{0}
\setcounter{corollary}{0} \setcounter{remark}{0}
\setcounter{conjecture}{0}

In view of (\ref{3.10}), our following conjecture implies that for any prime $p=x^2+7y^2$
with $x,y\in\Z^+$ we have
$$\sum_{k=0}^{p-1}\f{\bi{2k}k^2T_k(1,16)}{(-256)^k}\eq\l(\f{-1}p\r)(4x^2-2p)\pmod{p^2}.$$

\begin{conjecture}\label{Conj4.1} Let $p$ be an odd prime with $(\f p7)=1$.
Write $p=x^2+7y^2$ with $x,y\in\Z$ such that $x\eq1\ (\mo\ 4)$ if
$p\eq1\ (\mo\ 4)$, and $y\eq1\ (\mo\ 4)$ if $p\eq3\ (\mo\ 4)$. Then
$$\sum_{k=0}^{p-1}\f{\bi{2k}k^2}{256^k}u_k(1,16)\eq\begin{cases}0\ (\mo\ p^2)&\t{if}\ p\eq1\ (\mo\ 4),
\\\f13(\f2p)(\f p{7y}-4y)\ (\mo\ p^2)&\t{if}\ p\eq3\ (\mo\ 4);
\end{cases}$$
$$\sum_{k=0}^{p-1}\f{\bi{2k}k^2}{256^k}v_k(1,16)\eq\begin{cases}2(\f 2p)(2x-\f p{2x})\ (\mo\ p^2)&\t{if}\ p\eq1\ (\mo\ 4),
\\0\ (\mo\ p^2)&\t{if}\ p\eq3\ (\mo\ 4).
\end{cases}$$
When $p\eq1\ (\mo\ 4)$, we have
$$\sum_{k=0}^{p-1}\f{k\bi{2k}k^2}{16^k}u_k(1,16)\eq \sum_{k=0}^{p-1}\f{k\bi{2k}k^2}{256^k}u_k(1,16)
\eq\f{(\f2p)}{42}\l(x-\f p{2x}\r)\ (\mo\ p^2)$$
and
\begin{align*}&\sum_{k=0}^{p-1}(4k+3)\f{\bi{2k}k^2}{16^k}v_k(1,16)
\\\eq&3\sum_{k=0}^{p-1}(4k+1)\f{\bi{2k}k^2}{256^k}v_k(1,16)\eq6\l(\f 2p\r)x\ (\mo\ p^2).
\end{align*}
When $p\eq3\ (\mo\ 4)$, we can determine $y$ mod $p^2$ in the
following way:
$$\sum_{k=0}^{p-1}\f{k\bi{2k}k^2}{16^k}u_k(1,16)\eq\sum_{k=0}^{p-1}\f{k\bi{2k}k^2}{16^k}v_k(1,16)
\eq-\l(\f 2p\r)\f y2\pmod{p^2}$$
and
$$3\sum_{k=0}^{p-1}\f{k\bi{2k}k^2}{256^k}u_k(1,16)\eq\sum_{k=0}^{p-1}\f{k\bi{2k}k^2}{256^k}v_k(1,16)
\eq\l(\f 2p\r)\f y2\pmod{p^2}.$$
\end{conjecture}

Just like $\Q(\sqrt{-7})$, the imaginary quadratic field
$\Q(\sqrt{-11})$ also has class number one. Let $p$ be an odd prime.
Whenever $(\f p{11})=1$ we can write $p$ in the form $(x^2+11y^2)/4$ with $x,y\in\Z$.
We guess that
$$\sum_{k=0}^{p-1}\f{\bi{2k}k\bi{3k}kT_k(46,1)}{512^k}\eq\begin{cases} x^2-2p\ (\mo\ p^2)&\t{if}\ 4p=x^2+11y^2,
\\0\ (\mo\ p^2)&\t{if}\ (\f p{11})=-1.\end{cases}$$
To attack this we note that (\ref{3.12}) with $b=46$, $c=1$ and
$x=-27/512$ yields
$$\sum_{k=0}^{p-1}\f{\bi{2k}k\bi{3k}k}{512^k}T_k(46,1)
\eq\(\sum_{k=0}^{p-1}\f{\bi{2k}k\bi{3k}k}{(-64)^k}\al^k\)\times \sum_{k=0}^{p-1}\f{\bi{2k}k\bi{3k}k}{(-64)^k}\beta^k
\pmod{p^2},$$
where $\al=(1+\sqrt{33})/2$ and $\beta=(1-\sqrt{33})/2$. Observe that
$2\al^k=v_k(1,-8)+(\al-\beta)u_k(1,-8)$ and
$2\beta^k=v_k(1,-8)-(\al-\beta)u_k(1,-8).$ So we have
\begin{align*}&4\sum_{k=0}^{p-1}\f{\bi{2k}k\bi{3k}kT_k(46,1)}{512^k}
\\\eq&\(\sum_{k=0}^{p-1}\f{\bi{2k}k\bi{3k}k}{(-64)^k}v_k(1,-8)\)^2
-33\(\sum_{k=0}^{p-1}\f{\bi{2k}k\bi{3k}k}{(-64)^k}v_k(1,-8)\)^2\pmod{p^2}.
\end{align*}
This, together with the author's conjecture on
$\sum_{k=0}^{p-1}\bi{2k}k^2\bi{3k}k/64^k$ mod $p^2$ (cf. \cite[Conjecture 5.4]{Su2}) leads us to raise the following conjecture.

\begin{conjecture}\label{Conj4.2} Let $p>3$ be a prime. If $(\f p{11})=-1$,
then
$$\sum_{k=0}^{p-1}\f{\bi{2k}k\bi{3k}k}{(-64)^k}u_k(1,-8)\eq \sum_{k=0}^{p-1}\f{\bi{2k}k\bi{3k}k}{(-64)^k}v_k(1,-8)\eq0\pmod{p}.$$
When $(\f p{11})=1$, $p\eq1\ (\mo\ 3)$, and $4p=x^2+11y^2$
with $x\eq1\ (\mo\ 3)$, we have
\begin{align*}\sum_{k=0}^{p-1}\f{\bi{2k}k\bi{3k}k}{(-64)^k}u_k(1,-8)
\eq&0\ (\mo\ p^2),
\\\sum_{k=0}^{p-1}\f{k\bi{2k}k\bi{3k}k}{(-64)^k}u_k(1,-8)
\eq&\f{114}{11}\l(\f{2p}x-x\r)\ (\mo\ p^2),
\\\sum_{k=0}^{p-1}\f{k\bi{2k}k\bi{3k}k}{216^k}u_k(8,27)
\eq&\f{4}{99}\l(\f{2p}x-x\r)\ (\mo\ p^2),
\\\sum_{k=0}^{p-1}\f{\bi{2k}k\bi{3k}k}{(-64)^k}v_k(1,-8)
\eq&\sum_{k=0}^{p-1}\f{\bi{2k}k\bi{3k}k}{216^k}v_k(8,27)\eq 2\l(\f px-x\r)\pmod{p^2},
\end{align*}
and
\begin{align*}\sum_{k=0}^{p-1}(k+60)\f{\bi{2k}k\bi{3k}k}{(-64)^k}v_k(1,-8)\eq&-60x\pmod{p^2},
\\\sum_{k=0}^{p-1}(9k+2)\f{\bi{2k}k\bi{3k}k}{216^k}v_k(8,27)\eq&-2x\pmod{p^2}.
\end{align*}
When  $(\f p{11})=1$, $p\eq2\ (\mo\ 3)$, and $4p=x^2+11y^2$
with $y\eq1\ (\mo\ 3)$, we have
\begin{align*} &11\sum_{k=0}^{p-1}\f{\bi{2k}k\bi{3k}k}{(-64)^k}u_k(1,-8)
\\\eq&-3\sum_{k=0}^{p-1}\f{\bi{2k}k\bi{3k}k}{(-64)^k}v_k(1,-8)\eq\f 32\l(\f p{y}-11y\r)\ (\mo\ p^2),
\end{align*}
\begin{align*}\sum_{k=0}^{p-1}(2k-155)\f{\bi{2k}k\bi{3k}k}{(-64)^k}u_k(1,-8)\eq&\f{759}2y\pmod{p^2},
\\\sum_{k=0}^{p-1}(2k-243)\f{\bi{2k}k\bi{3k}k}{(-64)^k}v_k(1,-8)\eq&-\f{4359}2y\pmod{p^2},
\end{align*}
\begin{gather*}
\sum_{k=0}^{p-1}\f{\bi{2k}k\bi{3k}k}{216^k}u_k(8,27)\eq y-\f p{11y}\pmod{p^2},
\\\sum_{k=0}^{p-1}\f{k\bi{2k}k\bi{3k}k}{216^k}u_k(8,27)\eq\f18\sum_{k=0}^{p-1}\f{k\bi{2k}k\bi{3k}k}{216^k}v_k(8,27)\eq-\f y9\pmod{p^2}.
\end{gather*}
\end{conjecture}

Motivated by the author's investigation of
$\sum_{k=0}^{p-1}\bi{2k}k\bi{3k}kT_k(3,1)/27^k$ mod $p^2$ (with
$p>3$ a prime) and the congruence (\ref{3.12}), we pose the following
conjecture which involves the well-known Fibonacci numbers
$F_k=u_k(1,-1)\ (k\in\N)$ and Lucas numbers $L_k=v_k(1,-1)\
(k\in\N)$. Note that the imaginary quadratic field $\Q(\sqrt{-15})$
has class number 2.

\begin{conjecture}\label{Conj4.3} Let $p>5$ be a prime. If $p\eq1,4\ (\mo\
15)$ and $p=x^2+15y^2\ (x,y\in\Z)$ with $x\eq1\ (\mo\ 3)$, then
\begin{align*}\sum_{k=0}^{p-1}\f{k\bi{2k}k\bi{3k}k}{27^k}F_k\eq&\f 2{15}\l(\f px-2x\r)\pmod{p^2},
\\\sum_{k=0}^{p-1}\f{\bi{2k}k\bi{3k}k}{27^k}L_k\eq& 4x-\f px\pmod{p^2}
\end{align*}
and
$$\sum_{k=0}^{p-1}(3k+2)\f{\bi{2k}k\bi{3k}k}{27^k}L_k\eq 4x\pmod{p^2}.$$
If $p\eq2,8\ (\mo\ 15)$ and $p=3x^2+5y^2\ (x,y\in\Z)$ with $y\eq1\ (\mo\ 3)$, then
$$\sum_{k=0}^{p-1}\f{\bi{2k}k\bi{3k}k}{27^k}F_k\eq\f p{5y}-4y\pmod{p^2}$$
and
$$\sum_{k=0}^{p-1}\f{k\bi{2k}k\bi{3k}k}{27^k}F_k\eq\sum_{k=0}^{p-1}\f{k\bi{2k}k\bi{3k}k}{27^k}L_k\eq\f 43y\pmod{p^2}.$$
\end{conjecture}
\begin{remark}\label{Rem4.1} By \cite[Theorem 1.6]{Su8}, for any prime $p>3$ we have
$$\sum_{k=0}^{p-1}\f{\bi{2k}k\bi{3k}k}{27^k}F_k\eq0\ (\mo\ p^2)\ \ \t{if}\ \ p\eq1\ (\mo\ 3),$$
and
$$\sum_{k=0}^{p-1}\f{\bi{2k}k\bi{3k}k}{27^k}L_k\eq0\ (\mo\ p^2)\ \ \t{if}\ \ p\eq2\ (\mo\ 3).$$
\end{remark}

In fact, we have many other conjectures similar to Conjectures
4.1-4.3; for the sake of brevity we don't include them in this
paper.

\begin{conjecture}\label{Conj4.4} Let $p>3$ be a prime.

{\rm (i)} If $p\eq1,4\pmod{15}$ and $p=x^2+15y^2$ with $x,y\in\Z$,
then
$$P_{(p-1)/2}(7\sqrt{-15}\pm16\sqrt{-3})\eq\l(\f{-\sqrt{-15}}p\r)\l(\f x{15}\r)\l(2x-\f p{2x}\r)\pmod{p^2}.$$

{\rm (ii)} Suppose that $(\f p5)=(\f p{7})=1$ and write
$4p=x^2+35y^2$ with $x,y\in\Z$. If $p\eq1\ (\mo\ 3)$, then
$$\sum_{k=0}^{p-1}\f{\bi{2k}{k}\bi{3k}k}{3456^k}(64+27\sqrt{5}\pm \sqrt{-35})^k\eq\l(\f x3\r)\l(\f p{x}-x\r)\pmod{p^2}.$$
If $p\eq2\ (\mo\ 3)$, then
$$\sum_{k=0}^{p-1}\f{\bi{2k}{k}\bi{3k}k}{3456^k}(64+27\sqrt{5}\pm \sqrt{-35})^k\eq \pm\sqrt{-35}\l(\f y3\r)\l(y-\f p{35y}\r)\pmod{p^2}.$$

{\rm (iii)} If $(\f 2p)=(\f p3)=(\f p5)=1$ and $p=x^2+30y^2$ with
$x,y\in\Z$, then
$$\sum_{k=0}^{p-1}\f{\bi{2k}k\bi{3k}k}{2916^k}(54-35\sqrt2\pm\sqrt5)^k\eq\l(\f x3\r)\l(2x-\f p{2x}\r)\pmod{p^2}.$$

{\rm (iv)} If $(\f{-2}p)=(\f p3)=(\f p7)=1$, and $p=x^2+42y^2$ with
$x,y\in\Z$, then
$$\sum_{k=0}^{p-1}\f{\bi{2k}k\bi{3k}k}{13500^k}(250-99\sqrt6\pm 2\sqrt{14})^k\eq\l(\f x3\r)\l(2x-\f p{2x}\r)\pmod{p^2}.$$

{\rm (v)} If $(\f 2p)=(\f p3)=(\f p{13})=1$ and $p=x^2+78y^2$ with
$x,y\in\Z$, then
$$\sum_{k=0}^{p-1}\f{\bi{2k}k\bi{3k}k}{530604^k}(9826-6930\sqrt2\pm 5\sqrt{26})^k\eq\l(\f x3\r)\l(2x-\f p{2x}\r)\pmod{p^2}.$$

{\rm (vi)} If $(\f 2p)=(\f p3)=(\f p{17})=1$ and $p=x^2+102y^2$ with
$x,y\in\Z$, then
$$\sum_{k=0}^{p-1}\f{\bi{2k}k\bi{3k}k}{3881196^k}(71874-17420\sqrt{17}\pm 35\sqrt{2})^k\eq\l(\f x3\r)\l(2x-\f p{2x}\r)\pmod{p^2}.$$

{\rm (vii)} If $(\f {-1}p)=(\f p3)=(\f p{11})=1$ and $p=x^2+33y^2$
with $x,y\in\Z$, then
$$\sum_{k=0}^{p-1}\f{\bi{4k}{2k}\bi{2k}k}{(2^{12}3)^k}(96-5\sqrt{11}\pm 65\sqrt{3})^k\eq\l(\f x3\r)\l(\f p{2x}-2x\r)\pmod{p^2}.$$
\end{conjecture}

\begin{remark}\label{Rem4.2} Let $p\eq1,4\pmod{15}$ be a prime with $p=x^2+15y^2$
($x,y\in\Z$). Applying (\ref{3.10}) we see that
\begin{align*}&\sum_{k=0}^{p-1}(-1)^k\bi{2k}k^2T_k
\\\eq& P_{(p-1)/2}(7\sqrt{-15}+16\sqrt{-3}) P_{(p-1)/2}(7\sqrt{-15}-16\sqrt{-3})\pmod{p^2}.\end{align*}
Thus part (i) of Conjecture 4.4 implies that
$$\sum_{k=0}^{p-1}(-1)^k\bi{2k}k^2T_k\eq\l(2x-\f p{2x}\r)^2\eq4x^2-2p\pmod{p^2}.$$
We omit here similar comments on parts (ii)-(vii) of Conjecture 4.4.
We also have many other conjectures similar to Conjecture 4.4.
\end{remark}

\begin{conjecture}\label{Conj4.5} Let $p>5$ be a prime. Then
\begin{align*}&\l(\f{-1}p\r)\sum_{k=0}^{p-1}\f{\bi{2k}k^2T_{2k}(62,1)}{(-128^2)^k}
\eq\l(\f p3\r)\sum_{k=0}^{p-1}\f{\bi{2k}k^2T_{2k}(62,1)}{(-480^2)^k}
\\\eq&\begin{cases} 4x^2-2p\pmod{p^2}&\t{if}\ p\eq1,9\pmod{20}\ \&\ p=x^2+5y^2,
\\2x^2-2p\pmod{p^2}&\t{if}\ p\eq3,7\pmod{20}\ \&\ 2p=x^2+5y^2,
\\0\ (\mo\ p^2)&\t{if}\  p\eq 11,13,17,19\,(\mo\ 20).\end{cases}
\end{align*}
And
\begin{align*} \sum_{k=0}^{p-1}(340k+111)\f{\bi{2k}k^2T_{2k}(62,1)}{(-128^2)^k}\eq& 3p\l(\f {-1}p\r)\l(22+15\l(\f{p}{15}\r)\r)\pmod{p^2},
\\\sum_{k=0}^{p-1}(340k+59)\f{\bi{2k}k^2T_{2k}(62,1)}{(-480^2)^k}\eq& p\l(\f {-1}p\r)\l(51+8\l(\f{p}{15}\r)\r)\pmod{p^2}.
\end{align*}
\end{conjecture}

\begin{conjecture}\label{Conj4.6} Let $p>3$ be a prime.
 Then
 \begin{align*}&\sum_{k=0}^{p-1}\f{\bi{2k}kT_k^2}{4^k}\eq\l(\f p3\r)\sum_{k=0}^{p-1}\f{\bi{2k}kT_k^2(4,1)}{16^k}
\\\eq&\sum_{k=0}^{p-1}\f{\bi{2k}k^2T_k(10,1)}{(-64)^k}\eq
 \l(\f p3\r)\sum_{k=0}^{p-1}\f{\bi{2k}k^2T_{2k}(6,1)}{256^k}\eq\sum_{k=0}^{p-1}\f{\bi{2k}k^2T_{2k}(6,1)}{1024^k}
 \\\eq&\begin{cases}
 4x^2-2p\, (\mo\ p^2)&\t{if}\ p\eq1,7\,(\mo\ 24),\ p=x^2+6y^2,
 \\8x^2-2p\, (\mo\ p^2)&\t{if}\ p\eq5,11\, (\mo\ 24),\ p=2x^2+3y^2,
 \\0\, (\mo\ p^2)&\t{if}\ (\f {-6}p)=-1;
 \end{cases}\end{align*}
 and
\begin{align*}&\sum_{k=0}^{p-1}\f{\bi{2k}kT_k^2(6,1)}{192^k}
\\\eq&\begin{cases}
 (\f{-1}p)(4x^2-2p)\,(\mo\ p^2)&\t{if}\ p\eq1,7\, (\mo\ 24),\ p=x^2+6y^2,
 \\8x^2-2p\, (\mo\ p^2)&\t{if}\ p\eq5,11\, (\mo\ 24),\ p=2x^2+3y^2,
 \\0\, (\mo\ p^2)&\t{if}\ (\f {-6}p)=-1.
 \end{cases}\end{align*}
 Also,
 \begin{align*}\sum_{k=0}^{p-1}(3k+1)\f{\bi{2k}k^2T_k(10,1)}{(-64)^k}\eq&\f
 p4\l(3\l(\f p3\r)+1\r)\pmod{p^2},
 \\\sum_{k=0}^{p-1}(4k+1)\f{\bi{2k}kT_k^2(6,1)}{192^k}\eq& p\l(\f{-6}p\r)\l(4-3\l(\f 2p\r)\r)\pmod{p^2}.
 \end{align*}
\end{conjecture}

\begin{conjecture}\label{Conj4.7} Let $p>5$ be a prime.

 {\rm (i)} We have
 \begin{align*}&\sum_{k=0}^{p-1}\f{\bi{2k}kT_k^2(3,1)}{36^k}
\eq\l(\f2p\r)\sum_{k=0}^{p-1}\f{\bi{2k}k^2T_k(34,1)}{(-64)^k}
 \eq\sum_{k=0}^{p-1}\f{\bi{2k}k^2T_{2k}(18,1)}{4096^k}
  \\\eq&\begin{cases}
 4x^2-2p\, (\mo\ p^2)&\t{if}\ p\eq1,9,11,19\,(\mo\ 40), \, p=x^2+10y^2,
 \\8x^2-2p\,(\mo\ p^2)&\t{if}\ p\eq7,13,23,37\,(\mo\ 40),\, p=2x^2+5y^2,
 \\0\, (\mo\ p^2)&\t{if}\ (\f {-10}p)=-1.
 \end{cases}\end{align*}
Also,
\begin{align*} \sum_{k=0}^{p-1}(16k+5)\f{\bi{2k}kT_k(3,1)^2}{36^k}\eq&5p\pmod{p^2},
\\\sum_{k=0}^{p-1}(60k+23)\f{\bi{2k}k^2T_k(34,1)}{(-64)^k}\eq& p\l(8\l(\f2p\r)+15\l(\f{-1}p\r)\r)\pmod{p^2}.
\end{align*}
\end{conjecture}

\begin{conjecture}\label{Conj4.8} Let $p>7$ be a prime. Then
\begin{align*}&\sum_{k=0}^{p-1}\f{\bi{2k}k\bi{3k}kT_k(18,1)}{512^k}
\eq\l(\f{10}p\r)\sum_{k=0}^{p-1}\f{\bi{2k}k\bi{3k}kT_{3k}(6,1)}{(-512)^k}
\\\eq&
\begin{cases} x^2-2p\pmod{p^2}&\t{if}\ (\f p5)=(\f p7)=1\ \&\ 4p=x^2+35y^2,
\\ 2p-5x^2\pmod{p^2}&\t{if}\ (\f p5)=(\f p7)=-1\ \&\ 4p=5x^2+7y^2,
\\0\pmod{p^2}&\t{if}\ (\f p{35})=-1.\end{cases}
\end{align*}
And
\begin{align*}\sum_{k=0}^{p-1}(35k+9)\f{\bi{2k}k\bi{3k}kT_k(18,1)}{512^k}
\eq&\f{9p}2\l(7-5\l(\f p5\r)\r)\pmod{p^2},
\\\sum_{k=0}^{p-1}(35k+9)\f{\bi{2k}k^2T_{3k}(6,1)}{(-512)^k}\eq&\f {9p}{32}\l(\f 2p\r)\l(25+7\l(\f{p}7\r)\r)\pmod{p^2}.
\end{align*}
\end{conjecture}

\begin{conjecture}\label{Conj4.9} Let $p\not=2,29$ be a prime. When
$p\not=5,7$, we have
\begin{align*}&\sum_{k=0}^{p-1}\f{\bi{2k}k^2T_{2k}(19602,1)}{78400^{2k}}
\\\eq&\begin{cases} 4x^2-2p\pmod{p^2}&\t{if}\ (\f{-2}p)=(\f{29}p)=1\ \&\ p=x^2+58y^2,
\\8x^2-2p\pmod{p^2}&\t{if}\ (\f{-2}p)=(\f{29}p)=-1\ \&\ p=2x^2+29y^2,
\\0\pmod{p^2}&\t{if}\ (\f{-58}p)=-1.\end{cases}
\end{align*}
Provided $p\not=13$ we have
\begin{align*}&\sum_{k=0}^{p-1}\f{\bi{2k}k^2T_{2k}(19602,1)}{78416^{2k}}
\\\eq&\begin{cases} 4x^2-2p\pmod{p^2}&\t{if}\ (\f{-2}p)=(\f{29}p)=1\ \&\ p=x^2+58y^2,
\\2p-8x^2\pmod{p^2}&\t{if}\ (\f{-2}p)=(\f{29}p)=-1\ \&\ p=2x^2+29y^2,
\\0\pmod{p^2}&\t{if}\ (\f{-58}p)=-1.\end{cases}
\end{align*}
\end{conjecture}

\begin{conjecture}\label{Conj4.10} Let $p>5$ be a prime. Then
\begin{align*}&\l(\f{-6}p\r)\sum_{k=0}^{p-1}\f{\bi{2k}k\bi{3k}kT_{3k}(26,1)}{(-24)^{3k}}
\eq\l(\f{15}p\r)\sum_{k=0}^{p-1}\f{\bi{2k}k\bi{3k}kT_{3k}(62,1)}{(-240)^{3k}}
\\\eq&
\begin{cases} x^2-2p\pmod{p^2}&\t{if}\ (\f p7)=(\f p{13})=1\ \&\
4p=x^2+91y^2,
\\2p-7x^2\pmod{p^2}&\t{if}\ (\f p7)=(\f p{13})=-1\ \&\ 4p=7x^2+13y^2,
\\0\pmod{p^2}&\t{if}\ (\f p{91})=-1.\end{cases}
\end{align*}
And
\begin{align*}&\sum_{k=0}^{p-1}(819k+239)\f{\bi{2k}k\bi{3k}kT_{3k}(26,1)}{(-24)^{3k}}
\\\eq&\f {p}{32}\l(\f {-6}p\r)\l(949+6699\l(\f{p}7\r)\r)\pmod{p^2},
\end{align*}
\begin{align*}
&\sum_{k=0}^{p-1}(1638k+277)\f{\bi{2k}k\bi{3k}kT_{3k}(62,1)}{(-240)^{3k}}
\\\eq&\f p{40}\l(\f{-105}p\r)\l(8701+2379\l(\f p7\r)\r)\pmod{p^2}.
\end{align*}
\end{conjecture}
\begin{remark}\label{Rem4.3} Note that the imaginary quadratic field
$\Q(\sqrt{-d})$ has class number two for $d=5,6,10,15,35,58,91$.
\end{remark}

\begin{conjecture}\label{Conj4.11} Let $p>3$ be a prime. We have
\begin{align*}&\l(\f{-6}p\r)\sum_{k=0}^{p-1}\f{\bi{4k}{2k}\bi{2k}kT_{k}(110,1)}{(-96^2)^k}
\\\eq&\begin{cases} 4x^2-2p\pmod{p^2}&\t{if}\ (\f{-1}p)=(\f{p}3)=(\f p{7})=1,\ p=x^2+21y^2,
\\12x^2-2p\,(\mo\ p^2)&\t{if}\ (\f{-1}p)=(\f p7)=-1, (\f p3)=1,\, p=3x^2+7y^2,
\\2x^2-2p\,(\mo\ p^2)&\t{if}\ (\f {-1}p)=(\f p{3})=-1, (\f p{7})=1,\, 2p=x^2+21y^2,
\\6x^2-2p\,(\mo\ p^2)&\t{if}\ (\f {-1}p)=1, (\f p{3})=(\f p7)=-1,\, 2p=3x^2+7y^2,
\\0\,(\mo\ p^2)&\t{if}\ (\f{-21}p)=-1,\end{cases}
\end{align*}
and
$$\sum_{k=0}^{p-1}(28k+5)\f{\bi{4k}{2k}\bi{2k}kT_{k}(110,1)}{(-96^2)^{k}}
\eq \f p8\l(\f{-6}p\r)\l(33+7\l(\f p{7}\r)\r)\pmod{p^2}.$$
\end{conjecture}

\begin{conjecture}\label{Conj4.12} Let $p>3$ be a prime. Then
\begin{align*}&\sum_{k=0}^{p-1}\f{\bi{2k}k^2 T_{2k}(18,1)}{256^k}
\\\eq&\begin{cases} 4x^2-2p\,(\mo\ p^2)&\t{if}\ (\f 2p)=(\f p3)=(\f p5)=1,\ p=x^2+30y^2,
\\12x^2-2p\,(\mo\ p^2)&\t{if}\ (\f p3)=1, (\f 2p)=(\f p5)=-1,\, p=3x^2+10y^2,
\\2p-8x^2\,(\mo\ p^2)&\t{if}\ (\f 2p)=1, (\f p3)=(\f p5)=-1,\, p=2x^2+15y^2,
\\2p-6x^2\,(\mo\ p^2)&\t{if}\ (\f p5)=1, (\f 2p)=(\f p3)=-1,\, 2p=3x^2+10y^2,
\\0\,(\mo\ p^2)&\t{if}\ (\f{-30}p)=-1.\end{cases}
\end{align*}
And
\begin{align*}&\sum_{k=0}^{p-1}\f{\bi{2k}k^2 T_{2k}(30,1)}{256^k}
\\\eq&\begin{cases} 4x^2-2p\,(\mo\ p^2)&\t{if}\ (\f {-2}p)=(\f p3)=(\f p7)=1,\ p=x^2+42y^2,
\\12x^2-2p\,(\mo\ p^2)&\t{if}\ (\f {-2}p)=1, (\f p3)=(\f p7)=-1,\, p=3x^2+14y^2,
\\2p-8x^2\,(\mo\ p^2)&\t{if}\ (\f p7)=1, (\f {-2}p)=(\f p3)=-1,\, p=2x^2+21y^2,
\\2p-6x^2\,(\mo\ p^2)&\t{if}\ (\f p3)=1, (\f {-2}p)=(\f p7)=-1,\, 2p=3x^2+14y^2,
\\0\,(\mo\ p^2)&\t{if}\ (\f{-42}p)=-1.\end{cases}
\end{align*}
\end{conjecture}

\begin{conjecture}\label{Conj4.13} Let $p>3$ be a prime. When $p\not=13,17$,
we have
\begin{align*}&\sum_{k=0}^{p-1}\f{\bi{2k}k\bi{3k}kT_k(102,1)}{102^{3k}}
\\\eq&\begin{cases} 4x^2-2p\pmod{p^2}&\t{if}\ (\f 2p)=(\f p3)=(\f p{13})=1,\ p=x^2+78y^2,
\\2p-8x^2\,(\mo\ p^2)&\t{if}\ (\f 2p)=1, (\f p3)=(\f p{13})=-1,\, p=2x^2+39y^2,
\\12x^2-2p\,(\mo\ p^2)&\t{if}\ (\f p{13})=1, (\f 2p)=(\f p3)=-1,\, p=3x^2+26y^2,
\\2p-24x^2\,(\mo\ p^2)&\t{if}\ (\f p3)=1, (\f 2p)=(\f p{13})=-1,\, p=6x^2+13y^2,
\\0\,(\mo\ p^2)&\t{if}\ (\f{-78}p)=-1.\end{cases}
\end{align*}
Provided $p\not=11,17$, we have
\begin{align*}&\sum_{k=0}^{p-1}\f{\bi{2k}k\bi{3k}kT_k(198,1)}{198^{3k}}
\\\eq&\begin{cases} 4x^2-2p\pmod{p^2}&\t{if}\ (\f {2}p)=(\f p3)=(\f p{17})=1,\ p=x^2+102y^2,
\\2p-8x^2\,(\mo\ p^2)&\t{if}\ (\f p{17})=1,(\f {2}p)=(\f p{3})=-1,\, p=2x^2+51y^2,
\\12x^2-2p\,(\mo\ p^2)&\t{if}\ (\f p{3})=1, (\f {2}p)=(\f p{17})=-1,\, p=3x^2+34y^2,
\\2p-24x^2\,(\mo\ p^2)&\t{if}\ (\f {2}p)=1, (\f p3)=(\f p{17})=-1,\, p=6x^2+17y^2,
\\0\,(\mo\ p^2)&\t{if}\ (\f{-102}p)=-1.\end{cases}
\end{align*}
\end{conjecture}

\begin{conjecture}\label{Conj4.14} Let $p$ be an odd prime and let
$m$ belong to the set $\{2,3,6,10,18,30,102,198\}$. If $p\nmid m$, then
\begin{equation}\label{4.1}\sum_{k=0}^{p-1}\f{\bi{2k}k\bi{3k}kT_k(m,1)}{m^{3k}}
\eq\sum_{k=0}^{p-1}\f{\bi{2k}k^2T_{2k}(m,1)}{256^k}\pmod{p^2}.\end{equation}
If $m^2\not\eq-12\ (\mo\ p)$, then
\begin{equation}\label{4.2}\sum_{k=0}^{p-1}\f{\bi{2k}k^2T_{2k}(m,1)}{256^k}
\eq\l(\f{m^2+12}p\r)\sum_{k=0}^{p-1}\f{\bi{4k}{2k}\bi{2k}kT_k(m^2-2,1)}{(m^2+12)^{2k}}
\pmod{p^2}.\end{equation}
\end{conjecture}
\begin{remark}\label{Rem4.4} We note that (\ref{4.1}) holds mod $p$ for any integer
$m\not\eq0\ (\mo\ p)$, and (\ref{4.2}) holds mod $p$ for any $m\in\Z$ with
$m^2\not\eq-12\ (\mo\ p)$.
\end{remark}

\begin{conjecture}\label{Conj4.15} Let $p\not=2,5,19$ be a prime. We have
\begin{align*}&\sum_{k=0}^{p-1}\f{\bi{2k}k^2T_{2k}(5778,1)}{1216^{2k}}
\\\eq&\begin{cases} 4x^2-2p\pmod{p^2}&\t{if}\ (\f{2}p)=(\f{p}5)=(\f p{19})=1,\ p=x^2+190y^2,
\\8x^2-2p\,(\mo\ p^2)&\t{if}\ (\f{2}p)=1, (\f p5)=(\f p{19})=-1,\, p=2x^2+95y^2,
\\2p-20x^2\,(\mo\ p^2)&\t{if}\ (\f 2p)=(\f p{5})=-1, (\f p{19})=1,\, p=5x^2+38y^2,
\\2p-40x^2\,(\mo\ p^2)&\t{if}\ (\f 2p)=(\f p{19})=-1, (\f p{5})=1,\, p=10x^2+19y^2,
\\0\,(\mo\ p^2)&\t{if}\ (\f{-190}p)=-1,\end{cases}
\end{align*}
and
\begin{align*}&\sum_{k=0}^{p-1}(57720k+24893)\f{\bi{2k}k^2T_{2k}(5778,1)}{1216^{2k}}
\\&\ \ \ \eq p\l(11548+13345\l(\f p{95}\r)\r)\pmod{p^2}.
\end{align*}
Provided $p\not=17$ we have
$$\sum_{k=0}^{p-1}\f{\bi{2k}k^2T_{2k}(5778,1)}{439280^{2k}}
\eq\l(\f p5\r)\sum_{k=0}^{p-1}\f{\bi{2k}k^2T_{2k}(5778,1)}{1216^{2k}}\pmod{p^2},$$
\begin{align*}&\sum_{k=0}^{p-1}(57720k+3967)\f{\bi{2k}k^2T_{2k}(5778,1)}{439280^{2k}}
\\&\ \ \ \eq p\l(\f p{19}\r)\l(3983-16\l(\f p{95}\r)\r)\pmod{p^2}.
\end{align*}
\end{conjecture}

\begin{conjecture}\label{Conj4.16} Let $p>5$ be a prime. Then
\begin{align*}&\sum_{k=0}^{p-1}\f{\bi{2k}k^2T_{2k}(198,1)}{224^{2k}}\eq\l(\f p7\r)\sum_{k=0}^{p-1}\f{\bi{2k}k^2T_{2k}(322,1)}{48^{4k}}
\\\eq&\begin{cases} 4x^2-2p\,(\mo\ p^2)&\t{if}\ (\f 2p)=(\f p5)=(\f p7)=1,\ p=x^2+70y^2,
\\8x^2-2p\,(\mo\ p^2)&\t{if}\ (\f p7)=1, (\f 2p)=(\f p5)=-1,\, p=2x^2+35y^2,
\\2p-20x^2\,(\mo\ p^2)&\t{if}\ (\f p5)=1, (\f 2p)=(\f p7)=-1,\, p=5x^2+14y^2,
\\28x^2-2p\,(\mo\ p^2)&\t{if}\ (\f 2p)=1, (\f p5)=(\f p7)=-1,\, p=7x^2+10y^2,
\\0\,(\mo\ p^2)&\t{if}\ (\f{-70}p)=-1.\end{cases}
\end{align*}
Also,
\begin{align*}&\sum_{k=0}^{p-1}\f{\bi{2k}k^2T_{2k}(322,1)}{(-2^{10}3^4)^k}
\\\eq&\begin{cases} 4x^2-2p\,(\mo\ p^2)&\t{if}\ (\f {-1}p)=(\f p5)=(\f p{17})=1,\ p=x^2+85y^2,
\\2p-2x^2\,(\mo\ p^2)&\t{if}\ (\f {p}{17})=1, (\f {-1}p)=(\f p5)=-1,\, 2p=x^2+85y^2,
\\2p-20x^2\,(\mo\ p^2)&\t{if}\ (\f {-1}p)=1, (\f p5)=(\f p{17})=-1,\, p=5x^2+17y^2,
\\10x^2-2p\,(\mo\ p^2)&\t{if}\ (\f p{5})=1, (\f {-1}p)=(\f p{17})=-1,\, 2p=5x^2+17y^2,
\\0\,(\mo\ p^2)&\t{if}\ (\f{-85}p)=-1.\end{cases}
\end{align*}
And
\begin{align*}&\sum_{k=0}^{p-1}\f{\bi{2k}k^2T_{2k}(1298,1)}{24^{4k}}
\\\eq&\begin{cases} 4x^2-2p\,(\mo\ p^2)&\t{if}\ (\f {-2}p)=(\f p5)=(\f p{13})=1,\ p=x^2+130y^2,
\\8x^2-2p\,(\mo\ p^2)&\t{if}\ (\f {-2}p)=1, (\f p5)=(\f p{13})=-1,\, p=2x^2+65y^2,
\\2p-20x^2\,(\mo\ p^2)&\t{if}\ (\f p5)=1, (\f {-2}p)=(\f p{13})=-1,\, p=5x^2+26y^2,
\\2p-40x^2\,(\mo\ p^2)&\t{if}\ (\f p{13})=1, (\f {-2}p)=(\f p5)=-1,\, p=10x^2+13y^2,
\\0\,(\mo\ p^2)&\t{if}\ (\f{-130}p)=-1.\end{cases}
\end{align*}
\end{conjecture}
\begin{remark}\label{Rem4.5}. The imaginary quadratic field $\Q(\sqrt{-d})$ has
class number four for $d=21,30,42,70,78,85,102,130,190$.
\end{remark}

 \begin{conjecture}\label{Conj4.17} Let $p>3$ be a prime. Then
$$\sum_{k=0}^{p-1}\f{\bi{2k}k\bi{3k}kT_{3k}}{(-27)^k}\eq\begin{cases}(\f
p3)(4x^2-2p)\pmod{p^2}&\t{if}\ p=x^2+7y^2,
\\0\pmod{p^2}&\t{if}\ (\f p7)=-1;\end{cases}$$
$$\sum_{k=0}^{p-1}\f{\bi{2k}k\bi{3k}kT_{3k}(26,81)}{24^{3k}}\eq\begin{cases}(\f
6p)(x^2-2p)\pmod{p^2}&\t{if}\ 4p=x^2+11y^2,
\\0\pmod{p^2}&\t{if}\ (\f p{11})=-1;\end{cases}$$
 $$\sum_{k=0}^{p-1}\f{\bi{2k}k\bi{3k}kT_{3k}(10,1)}{24^{3k}}
 \eq\begin{cases}
 (\f 6p)(x^2-2p)\pmod{p^2}&\t{if}\ 4p=x^2+19y^2,
 \\0\pmod{p^2}&\t{if}\ (\f p{19})=-1.\end{cases}$$
If $p\not=13$, then
$$\sum_{k=0}^{p-1}\f{\bi{2k}k\bi{3k}kT_{3k}(106,1)}{312^{3k}}
\eq\begin{cases}
 (\f {78}p)(x^2-2p)\pmod{p^2}&\t{if}\ 4p=x^2+43y^2,
 \\0\pmod{p^2}&\t{if}\ (\f p{43})=-1.\end{cases}$$
If $p\not=73$, then
\begin{align*}&\sum_{k=0}^{p-1}\f{\bi{2k}k\bi{3k}kT_{3k}(586,1)}{1752^{3k}}
\\\eq&\begin{cases}
 (\f {438}p)(x^2-2p)\ (\mo\ p^2)&\t{if}\  (\f p{67})=1\ \&\ 4p=x^2+67y^2,
 \\0\ (\mo\ p^2)&\t{if}\ (\f p{67})=-1.\end{cases}
\end{align*}
If $p\not=8893$, then
\begin{align*}&\sum_{k=0}^{p-1}\f{\bi{2k}k\bi{3k}kT_{3k}(71146,1)}{213432^{3k}}
\\\eq&\begin{cases}
 (\f {53358}p)(x^2-2p)\pmod{p^2}&\t{if}\ (\f p{163})=1\ \&\ 4p=x^2+163y^2,
 \\0\pmod{p^2}&\t{if}\ (\f p{163})=-1.\end{cases}
 \end{align*}
 Also,
 \begin{align*}&\sum_{k=0}^{p-1}\f{\bi{2k}k\bi{3k}kT_{3k}(2,-1)}{(-3456)^k}\eq
 \l(\f{-1}p\r)\sum_{k=0}^{p-1}\f{\bi{2k}k\bi{3k}kT_{3k}(2,9)}{24^{3k}}
 \\\eq&\begin{cases}
 (\f {2}p)(x^2-2p)\pmod{p^2}&\t{if}\ p\eq1\pmod3\ \&\ 4p=x^2+27y^2,
 \\0\pmod{p^2}&\t{if}\ p\eq2\pmod3.\end{cases}
 \end{align*}
And
\begin{align*}&\sum_{k=0}^{p-1}(15k+2)\f{\bi{2k}k\bi{3k}kT_{3k}(2,-1)}{(-3456)^k}
\\\eq&\begin{cases} 2p(\f2p)\pmod{p^2}&\t{if}\ 3\mid p-1\ \t{and}\ 2\ \t{is a cubis residue mod}\ p,
\\0\pmod{p}&\t{otherwise}.\end{cases}
\end{align*}
 \end{conjecture}
 \begin{remark}\label{Rem4.6} The imaginary quadratic field $\Q(\sqrt{-d})$ has class number one for $d=7,11,19,43,67,163$.
 We observe that if $p>3$ is a prime and $m$ is an integer with $m(3m+8)\not\eq0\ (\mo\ p)$ then
 \begin{align*}&\sum_{k=0}^{p-1}\f{\bi{2k}k\bi{3k}kT_{3k}(m+2,1)}{(3m)^{3k}}
 \\\eq&\l(\f{-m(3m+8)}p\r)\sum_{k=0}^{p-1}\f{\bi{2k}k\bi{3k}k\bi{6k}{3k}}{(-9m-24)^{3k}}\ \pmod{p}.
 \end{align*}
 \end{remark}

\begin{conjecture}\label{Conj4.18}
Let $p$ be an odd prime.

{\rm (i)} When $p>5$ we have
\begin{align*}&\sum_{k=0}^{p-1}\l(\f{T_k(38,21^2)}{(-16)^k}\r)^3\eq\l(\f{-5}p\r)\sum_{k=0}^{p-1}\l(\f{T_k(38,21^2)}{20^k}\r)^3
\\\eq&\begin{cases} 4x^2-2p\,(\mo\ p^2)&\t{if}\ (\f 2p)=(\f p3)=(\f p5)=1,\, p=x^2+30y^2,
\\12x^2-2p\,(\mo\ p^2)&\t{if}\ (\f p3)=1, (\f 2p)=(\f p5)=-1,\, p=3x^2+10y^2,
\\2p-8x^2\,(\mo\ p^2)&\t{if}\ (\f 2p)=1, (\f p3)=(\f p5)=-1,\, p=2x^2+15y^2,
\\20x^2-2p\,(\mo\ p^2)&\t{if}\ (\f p5)=1, (\f 2p)=(\f p3)=-1,\, p=5x^2+6y^2,
\\p \da_{p,7}\,(\mo\ p^2)&\t{if}\ (\f{-30}p)=-1.\end{cases}
\end{align*}
Also,
$$\sum_{k=0}^{p-1}(28k+15)\f{T_k^3(38,21^2)}{(-16)^{3k}}\eq \f p7\l(124-19\l(\f p3\r)\r)\pmod {p^2}.$$
If $p\not=7$, then
$$\sum_{k=0}^{p-1}\f{\bi{2k}kT_k^2(4,9)}{28^{2k}}\eq\sum_{k=0}^{p-1}\l(\f{T_k(38,21^2)}{(-16)^k}\r)^3\pmod{p^2}$$
and
$$\sum_{k=0}^{p-1}\f{24k+5}{28^{2k}}\bi{2k}kT_k^2(4,9)\eq p\l(\f{-6}p\r)\l(4+\l(\f 2p\r)\r)\pmod{p^2}.$$

{\rm (ii)} When $p\not=7$ we have
\begin{align*}&\sum_{k=0}^{p-1}\l(\f{T_k(110,57^2)}{32^k}\r)^3\eq\l(\f{-14}p\r)\sum_{k=0}^{p-1}\l(\f{T_k(110,57^2)}{(-28)^k}\r)^3
\\\eq&\begin{cases}4x^2-2p\,(\mo\ p^2)&\t{if}\ (\f{-2}p)=(\f p3)=(\f p{7})=1,\, p=x^2+42y^2,
\\8x^2-2p\,(\mo\ p^2)&\t{if}\ (\f{p}7)=1,(\f {-2}p)=(\f p{3})=-1,\, p=2x^2+21y^2,
\\12x^2-2p\,(\mo\ p^2)&\t{if}\ (\f{-2}{p})=1,(\f p3)=(\f p7)=-1,\, p=3x^2+14y^2,
\\24x^2-2p\,(\mo\ p^2)&\t{if}\ (\f{p}3)=1,(\f {-2}p)=(\f p{7})=-1,\,p=6x^2+7y^2,
\\p\da_{p,19}\,(\mo\ p^2)&\t{if}\ (\f{-42}p)=-1.\end{cases}
\end{align*}
Also,
$$\sum_{k=0}^{p-1}(684k+329)\f{T_k^3(110,57^2)}{2^{15k}}\eq \f p{19}\l(5160\l(\f{-2}p\r)+1091\r)\pmod{p^2}.$$
If $p>3$ and $p\not=11,19$, then
$$\sum_{k=0}^{p-1}\f{\bi{2k}kT_k^2(5,1)}{22^{2k}}\eq\sum_{k=0}^{p-1}\l(\f{T_k(110,57^2)}{32^k}\r)^3\pmod{p^2}$$
and
$$\sum_{k=0}^{p-1}\f{560k+71}{22^{2k}}\bi{2k}kT_k^2(5,1)\eq \f p3\l(280\l(\f p7\r)-67\r)\pmod{p^2}.$$
\end{conjecture}

\begin{conjecture}\label{Conj4.19} Let $p>3$ be a prime.
When $p\not=5$, we have
\begin{align*}&\sum_{k=0}^{p-1}\f{\bi{2k}kT_k^2(9,12)}{900^k}
\\\eq&\begin{cases}4x^2-2p\, (\mo\ p^2)&\t{if}\ (\f{-1}p)=(\f p3)=(\f p{11})=1,\ p=x^2+33y^2,
\\2x^2-2p\, (\mo\ p^2)&\t{if}\ (\f{-1}p)=1, (\f p3)=(\f p{11})=-1,\, 2p=x^2+33y^2,
\\12x^2-2p\, (\mo\ p^2)&\t{if}\ (\f{p}{11})=1, (\f {-1}p)=(\f p{3})=-1,\, p=3x^2+11y^2,
\\6x^2-2p\, (\mo\ p^2)&\t{if}\ (\f{p}{3})=1, (\f {-1}p)=(\f p{11})=-1,\, 2p=3x^2+11y^2,
\\0\, (\mo\ p^2)&\t{if}\ (\f{-33}p)=-1.\end{cases}
\end{align*}
Provided $p\not=29$, we have
\begin{align*}&\sum_{k=0}^{p-1}\f{\bi{2k}kT_k^2(171,-171)}{(-5177196)^k}
\\\eq&\begin{cases}4x^2-2p\,(\mo\ p^2)&\t{if}\ (\f{-1}p)=(\f p7)=(\f p{19})=1,\ p=x^2+133y^2,
\\2x^2-2p\,(\mo\ p^2)&\t{if}\ (\f p{7})=1, (\f{-1}p)=(\f p{19})=-1,\, 2p=x^2+133y^2,
\\2p-28x^2\,(\mo\ p^2)&\t{if}\ (\f{p}{19})=1, (\f{-1}p)=(\f p{7})=-1,\, p=7x^2+19y^2,
\\2p-14x^2\,(\mo\ p^2)&\t{if}\ (\f {-1}p)=1, (\f{p}7)=(\f p{19})=-1,\, 2p=7x^2+19y^2,
\\0\,(\mo\ p^2)&\t{if}\ (\f{-133}p)=-1.\end{cases}
\end{align*}
\end{conjecture}

\begin{conjecture}\label{Conj4.20} Let $p$ be an odd prime. When $p\not=23$,
we have
\begin{align*}&\sum_{k=0}^{p-1}\f{\bi{2k}kT_k^2(7,1)}{46^{2k}}
\\\eq&\begin{cases} 4x^2-2p\,(\mo\ p^2)&\t{if}\ (\f 2p)=(\f p5)=(\f p7)=1,\ p=x^2+70y^2,
\\8x^2-2p\,(\mo\ p^2)&\t{if}\ (\f p7)=1,\, (\f 2p)=(\f p5)=-1,\, p=2x^2+35y^2,
\\20x^2-2p\,(\mo\ p^2)&\t{if}\ (\f p5)=1,\, (\f 2p)=(\f p7)=-1,\, p=5x^2+14y^2,
\\2p-28x^2\,(\mo\ p^2)&\t{if}\ (\f 2p)=1,\, (\f p5)=(\f p7)=-1,\, p=7x^2+10y^2,
\\0\,(\mo\ p^2)&\t{if}\ (\f{-70}p)=-1.\end{cases}
\end{align*}
Provided $p\not=3,7,11,17,31$, we have
\begin{align*}&\sum_{k=0}^{p-1}\f{\bi{2k}kT_k^2(73,576)}{434^{2k}}
\\\eq&\begin{cases} 4x^2-2p\,(\mo\ p^2)&\t{if}\ (\f {2}p)=(\f p3)=(\f p{17})=1,\ p=x^2+102y^2,
\\8x^2-2p\,(\mo\ p^2)&\t{if}\ (\f p{17})=1, (\f {2}p)=(\f p{3})=-1,\ p=2x^2+51y^2,
\\12x^2-2p\,(\mo\ p^2)&\t{if}\ (\f p{3})=1, (\f {2}p)=(\f p{17})=-1,\ p=3x^2+34y^2,
\\24x^2-2p\,(\mo\ p^2)&\t{if}\ (\f {2}p)=1, (\f p3)=(\f p{17})=-1,\ p=6x^2+17y^2,
\\0\,(\mo\ p^2)&\t{if}\ (\f{-102}p)=-1.\end{cases}
\end{align*}
\end{conjecture}

\begin{conjecture}\label{Conj4.21} Let $p>3$ be a prime. If $p\not=23$, then
\begin{align*}&\sum_{k=0}^{p-1}\f{\bi{2k}kT_k^2(23,7^4)}{46^{2k}}
\\\eq&\begin{cases} 4x^2-2p\pmod{p^2}&\t{if}\ (\f 2p)=(\f p3)=(\f p{13})=1,\ p=x^2+78y^2,
\\8x^2-2p\,(\mo\ p^2)&\t{if}\ (\f 2p)=1, (\f p3)=(\f p{13})=-1,\, p=2x^2+39y^2,
\\2p-12x^2\,(\mo\ p^2)&\t{if}\ (\f p{13})=1, (\f 2p)=(\f p3)=-1,\, p=3x^2+26y^2,
\\2p-24x^2\,(\mo\ p^2)&\t{if}\ (\f p3)=1, (\f 2p)=(\f p{13})=-1,\, p=6x^2+13y^2,
\\p\da_{p,7}\,(\mo\ p^2)&\t{if}\ (\f{-78}p)=-1,\end{cases}
\end{align*}
where $\delta_{m,n}$ takes $1$ or $0$ according as $m=n$ or not.
If $p\not=5$, then
\begin{align*}&\sum_{k=0}^{p-1}\f{T_k^3(1298,651^2)}{(-100)^{3k}}
\\\eq&\begin{cases} 4x^2-2p\pmod{p^2}&\t{if}\ (\f 2p)=(\f p3)=(\f p{13})=1,\ p=x^2+78y^2,
\\2p-8x^2\,(\mo\ p^2)&\t{if}\ (\f 2p)=1, (\f p3)=(\f p{13})=-1,\, p=2x^2+39y^2,
\\12x^2-2p\,(\mo\ p^2)&\t{if}\ (\f p{13})=1,(\f 2p)=(\f p3)=-1,\, p=3x^2+26y^2,
\\2p-24x^2\,(\mo\ p^2)&\t{if}\ (\f p3)=1, (\f 2p)=(\f p{13})=-1,\, p=6x^2+13y^2,
\\p(\da_{p,7}+\da_{p,31})\pmod{p^2}&\t{if}\ (\f{-78}p)=-1.\end{cases}
\end{align*}
\end{conjecture}
\begin{conjecture}\label{Conj4.22} Let $p\not=2,7,11$ be a prime. Then
\begin{align*}&\sum_{k=0}^{p-1}\f{\bi{2k}kT_k^2(7,81)}{14^{2k}}\eq\l(\f p{11}\r)\sum_{k=0}^{p-1}\f{\bi{2k}kT_k^2(20,1)}{28^{2k}}
\\\eq&\begin{cases} 4x^2-2p\ (\mo\ p^2)&\t{if}\ (\f{-11}p)=(\f{2}p)=1\ \&\ p=x^2+22y^2,
\\8x^2-2p\ (\mo\ p^2)&\t{if}\ (\f{-11}p)=(\f{2}p)=-1\ \&\ p=2x^2+11y^2,
\\0\ (\mo\ p^2)&\t{if}\ (\f{-11}p)=-(\f{2}p).\end{cases}
\end{align*}
\end{conjecture}

\begin{conjecture}\label{Conj4.23} Let $p\not=2,7$ be  a prime. Then
$$\sum_{k=0}^{p-1}\f{\bi{2k}kT_k^2(6,2)}{450^k}\eq\begin{cases}4x^2-2p\ (\mo\ p^2)&\t{if}\ (\f p7)=1\ \&\ p=x^2+7y^2,
\\0\ (\mo\ p^2)&\t{if}\ (\f p7)=-1,\end{cases}$$
and
$$\sum_{k=0}^{p-1}\f{221k+28}{450^k}\bi{2k}kT_k^2(6,2)\eq\f{4p}7\l(72\l(\f{-1}p\r)-23\r)\pmod{p^2}.$$
\end{conjecture}

\begin{conjecture}\label{Conj4.24} Let $p>3$ be a prime. Then
\begin{align*}&\sum_{k=0}^{p-1}\f{\bi{2k}kT_k^2(7,12)}{4^k}\eq\sum_{k=0}^{p-1}\f{\bi{2k}kT_{2k}^2(3,3)}{36^k}
\\\eq&\begin{cases}4x^2-2p\ (\mo\ p^2)&\t{if}\ p\eq1\ (\mo\ 12)\ \&\ p=x^2+9y^2,
\\4xy\ (\mo\ p^2)&\t{if}\ p\eq5\ (\mo\ 12)\ \&\ p=x^2+y^2\ (3\mid x-y),
\\0\ (\mo\ p^2)&\t{if}\ p\eq3\ (\mo\ 4).\end{cases}
\end{align*}
\end{conjecture}

\begin{conjecture}\label{Conj4.25} Let $p$ be an odd prime. Then
\begin{align*}&\sum_{k=0}^{p-1}\f{\bi{2k}kT_k^2(19,-20)}{22^{2k}}\eq\sum_{k=0}^{p-1}\f{\bi{2k}kT_{2k}^2(9,20)}{4^k}
\\\eq&\begin{cases}4x^2-2p\ (\mo\ p^2)&\t{if}\ (\f{-1}p)=(\f p5)=1\ \&\ p=x^2+y^2\ (5\nmid x),
\\4xy\ (\mo\ p^2)&\t{if}\ (\f{-1}p)=-(\f p5)=1\ \&\ p=x^2+y^2\ (5\mid x-y),
\\0\ (\mo\ p^2)&\t{if}\ p\eq3\ (\mo\ 4)\ \&\ p\not=11.\end{cases}
\end{align*}
\end{conjecture}

We have many conjectures similar to Conjectures 4.18-4.25. For example, we find that
$\sum_{k=0}^{p-1}\bi{2k}kT_k^2(b,c)/m^k$ mod $p^2$ is related to the representation $p=x^2+dy^2$
if $(b,c,m;d)$ is among
\begin{gather*}(5,4,4;10),\ (3,-4,36;13),\ (5,4,14^2;30),
\\ (7,1,14^2;30),\ (7,28,14^2;21),\ (11,49,22^2;42).\end{gather*}

Though we will not list many other conjectures similar to
Conjectures 4.4-4.25, the above conjectures should convince the
reader that our conjectural series for $1/\pi$ in the next section
are indeed reasonable in view of the corresponding congruences.

\section{Dualities and new series for $1/\pi$}
\setcounter{lemma}{0} \setcounter{theorem}{0}
\setcounter{corollary}{0} \setcounter{remark}{0}
 \setcounter{conjecture}{0}

As mentioned in Section 1, for $b>0$ and $c>0$ the main term of
$T_n(b,c)$ as $n\to+\infty$ is
$$f_n(b,c):=\f{(b+2\sqrt{c})^{n+1/2}}{2\root{4}\of c\sqrt{n\pi}}.$$
Here we formulate a further refinement of this.

\begin{conjecture}\label{Conj5.1} For any positive real numbers $b$ and $c$,
we have
$$T_n(b,c)=f_n(b,c)\l(1+\f{b-4\sqrt{c}}{16n\sqrt{c}}+O\l(\f1{n^2}\r)\r)$$
as $n\to+\infty$. If $c>0$ and $b=4\sqrt{c}$, then
$$\f{T_n(b,c)}{\sqrt c^n}=T_n(4,1)=\f{3\times 6^n}{\sqrt{6n\pi}}\l(1+\f1{8n^2}+\f{15}{64n^3}+\f{21}{32n^4}+O\l(\f1{n^5}\r)\r).$$
If $c<0$ and $b\in\R$ then
$$\lim_{n\to\infty}\root{n}\of{|T_n(b,c)|}=\sqrt{b^2-4c}.$$
\end{conjecture}
\medskip

Let $p$ be an odd prime. Z.-H. Sun \cite{S1} proved the congruence
\begin{equation}\label{5.1}\sum_{k=0}^{p-1}\f{\bi{2k}k^2}{16^k}x^k\eq\l(\f{-1}p\r) \sum_{k=0}^{p-1}\f{\bi{2k}k^2}{16^k}(1-x)^k\pmod{p^2}\end{equation}
via Legendre polynomials; in fact this follows from the well-known
identity $P_n(-x)=(-1)^nP_n(x)$ with $n=(p-1)/2$. In \cite{Su8} the
author managed to show the following congruences via the Zeilberger
algorithm:
\begin{align}\label{5.2}\sum_{k=0}^{p-1}\f{\bi{2k}k\bi{4k}{2k}}{64^k}x^k
\eq&\l(\f{-2}p\r) \sum_{k=0}^{p-1}\f{\bi{2k}k\bi{4k}{2k}}{64^k}(1-x)^k\pmod{p^2},
\\\label{5.3}\sum_{k=0}^{p-1}\f{\bi{2k}k\bi{3k}{k}}{27^k}x^k
\eq&\l(\f{p}3\r)\sum_{k=0}^{p-1}\f{\bi{2k}k\bi{3k}{k}}{27^k}(1-x)^k\pmod{p^2}\ \ (p\not=3),
\\\label{5.4}\sum_{k=0}^{p-1}\f{\bi{3k}k\bi{6k}{3k}}{432^k}x^k
\eq&\l(\f{-1}p\r) \sum_{k=0}^{p-1}\f{\bi{3k}k\bi{6k}{3k}}{432^k}(1-x)^k\pmod{p^2}\ \ (p\not=3).
\end{align}

Our following result on dualities was motivated by
(38)-(41).

\begin{theorem}\label{Thm5.1} Let $p$ be an odd prime and let $b,c$ and
$m\not\eq0\pmod p$ be rational $p$-adic integers. Then
\begin{align}\label{5.5}\sum_{k=0}^{p-1}\f{\bi{2k}k^2}{(16m)^k}T_k(b,c)
\eq&\l(\f{-1}p\r)\sum_{k=0}^{p-1}\f{\bi{2k}k^2}{(16m)^k}T_k(m-b,c)\pmod{p^2},
\\\label{5.6}\sum_{k=0}^{p-1}\f{\bi{2k}k\bi{4k}{2k}}{(64m)^k}T_k(b,c)
\eq&\l(\f{-2}p\r)\sum_{k=0}^{p-1}\f{\bi{2k}k\bi{4k}{2k}}{(64m)^k}T_k(m-b,c)\pmod{p^2},
\end{align}
Provided $p>3$ we also have
\begin{align}\label{5.7}\sum_{k=0}^{p-1}\f{\bi{2k}k\bi{3k}k}{(27m)^k}T_k(b,c)
\eq&\l(\f{p}3\r)\sum_{k=0}^{p-1}\f{\bi{2k}k\bi{3k}k}{(27m)^k}T_k(m-b,c)\pmod{p^2},
\\\label{5.8}\sum_{k=0}^{p-1}\f{\bi{3k}k\bi{6k}{3k}}{(432m)^k}T_k(b,c)
\eq&\l(\f{-1}p\r)\sum_{k=0}^{p-1}\f{\bi{3k}k\bi{6k}{3k}}{(432m)^k}T_k(m-b,c)\pmod{p^2}.
\end{align}
\end{theorem}
\Proof.  Since the proofs of (42)-(45) are very similar, we
just show (43) in detail.

For $d=0,\ldots,p-1$, by taking differentiations of both sides (39)
$d$ times we get
$$\sum_{k=0}^{p-1}\f{\bi{2k}k\bi{4k}{2k}}{64^k}\bi kdx^{k-d}
\eq\l(\f{-2}p\r)\sum_{k=0}^{p-1}\f{\bi{2k}k\bi{4k}{2k}}{64^k}(-1)^d\bi kd(1-x)^{k-d}\pmod{p^2}.$$
In view of this, we have
\begin{align*}&\sum_{k=0}^{p-1}\f{\bi{2k}k\bi{4k}{2k}}{(64m)^k}T_k(b,c)
\\=&\sum_{k=0}^{p-1}\f{\bi{2k}k\bi{4k}{2k}}{(64m)^k}\sum_{j=0}^{\lfloor k/2\rfloor}\bi k{2j}\bi{2j}jb^{k-2j}c^j
\\=&\sum_{j=0}^{p-1}\bi{2j}j\f{c^j}{m^{2j}}\sum_{k=0}^{p-1}\f{\bi{2k}k\bi{4k}{2k}}{64^k}\bi k{2j}\l(\f bm\r)^{k-2j}
\\\eq&\sum_{j=0}^{p-1}\bi{2j}j\f{c^j}{m^{2j}}\l(\f{-2}p\r)
\sum_{k=0}^{p-1}\f{\bi{2k}k\bi{4k}{2k}}{64^k}\bi k{2j}\l(1-\f bm\r)^{k-2j}
\\=&\l(\f{-2}p\r)\sum_{k=0}^{p-1}\f{\bi{2k}k\bi{4k}{2k}}{(64m)^k}\sum_{j=0}^{\lfloor k/2\rfloor}\bi k{2j}\bi{2j}j(m-b)^{k-2j}c^j
\\=&\l(\f{-2}p\r)\sum_{k=0}^{p-1}\f{\bi{2k}k\bi{4k}{2k}}{(64m)^k}T_k(m-b,c)\pmod{p^2}.
\end{align*}

The proof of Theorem 5.1 is now complete. \qed

\medskip

\noindent{\it Example} 1. Let $p$ be an odd prime. By (42) we
have
$$\l(\f{-1}p\r)\sum_{k=0}^{p-1}\f{\bi{2k}k^2T_k(5,4)}{16^k}\eq\sum_{k=0}^{p-1}\f{\bi{2k}k^2T_k(-4,4)}{16^k}
=\sum_{k=0}^{p-1}\f{\bi{2k}k^3}{(-8)^k}\pmod{p^2}.$$ The author \cite{Su4}
conjectured that
$$\sum_{k=0}^{p-1}\f{\bi{2k}k^3}{(-8)^k}\eq\begin{cases} 4x^2-2p\ (\mo\ p^2)&\t{if}\ p=x^2+y^2\ (2\nmid x),
\\0\ (\mo\ p^2)&\t{if}\ p\eq3\ (\mo\ 4),\end{cases}$$
and this was recently confirmed by Z.-H. Sun \cite{S2}. When $p>3$, by
(44) we have
\begin{align*}&\l(\f p3\r)\sum_{k=0}^{p-1}\f{\bi{2k}k\bi{3k}kT_k(3,1)}{27^k}
\\\eq&\sum_{k=0}^{p-1}\f{\bi{2k}k\bi{3k}kT_k(-2,1)}{27^k}=\sum_{k=0}^{p-1}\f{\bi{2k}k^2\bi{3k}k}{(-27)^k}\pmod{p^2};
\end{align*}
the reader may consult \cite[Conjecture 5.6]{Su4} for
$\sum_{k=0}^{p-1}\bi{2k}k^2\bi{3k}k/(-27)^k$ mod $p^2$.
\medskip

Based on our investigations of congruences on sums of central
binomial coefficients and central trinomial coefficients, and the
author's philosophy about series for $1/\pi$ stated in \cite{Su7}, we
raise 61 conjectural series for $1/\pi$ of the following seven new
types with $a,b,c,d,m$ integers and $abcd(b^2-4c)m$ nonzero.

\ \ {\tt Type I}.
$\sum_{k=0}^\infty(a+dk)\bi{2k}k^2T_k(b,c)/m^k$.

\ \ {\tt Type II}.
$\sum_{k=0}^\infty(a+dk)\bi{2k}k\bi{3k}kT_k(b,c)/m^k$.

\ \ {\tt Type III}.
$\sum_{k=0}^\infty(a+dk)\bi{4k}{2k}\bi{2k}kT_k(b,c)/m^k$.

\ \ {\tt Type IV}.
$\sum_{k=0}^\infty(a+dk)\bi{2k}{k}^2T_{2k}(b,c)/m^k$.

\ \ {\tt Type V}.
$\sum_{k=0}^\infty(a+dk)\bi{2k}{k}\bi{3k}kT_{3k}(b,c)/m^k$.

\ \ {\tt Type VI}.
$\sum_{k=0}^\infty(a+dk)T_{k}(b,c)^3/m^k.$

\ \ {\tt Type VII}.
$\sum_{k=0}^\infty(a+dk)\bi{2k}kT_{k}(b,c)^2/m^k.$
\medskip

Recall that a series $\sum_{k=0}^\infty a_k$ is said to converge at
a geometric rate with ratio $r$ if
$\lim_{k\to+\infty}a_{k+1}/a_k=r\in(0,1)$. All the series in
Conjectures I-VII below converge at geometrical rates,
and they were found by the author in 2011 except that (IV19)-(IV21) were discovered in 2012.

\medskip

\noindent {\bf Conjecture I}. {\rm  We have the following identities:
\begin{align*}\sum_{k=0}^\infty\f{30k+7}{(-256)^k}\bi{2k}k^2T_k(1,16)=&\f{24}{\pi},\tag{I1}
\\\sum_{k=0}^\infty\f{30k+7}{(-1024)^k}\bi{2k}k^2T_k(34,1)=&\f{12}{\pi},\tag {I2}
\end{align*}
\begin{align*}
\\\sum_{k=0}^\infty\f{30k-1}{4096^k}\bi{2k}k^2T_k(194,1)=&\f{80}{\pi},\tag {I3}
\\\sum_{k=0}^\infty\f{42k+5}{4096^k}\bi{2k}k^2T_k(62,1)=&\f{16\sqrt3}{\pi}.\tag {I4}
\end{align*}}
\begin{remark}\label{Rem5.1} (I1) was the first identity for $1/\pi$ involving
generalized central trinomial coefficients; it was discovered by the author on
Jan. 2, 2011. Different from classical Ramanujan-type series for
$1/\pi$ (cf. N. D. Baruah and B. C. Berndt \cite{BB},
 and Berndt \cite[pp. 353-354]{Be})
and their known generalizations (see, e.g., S. Cooper \cite{C}), the two
coefficients in the linear part $30k-1$ of (I3) have {\it different
signs}, and also its corresponding $p$-adic congruence (with $p>3$ a
prime) involves {\it two} Legendre symbols:
$$\sum_{k=0}^{p-1}(30k-1)\f{\bi{2k}k^2T_k(194,1)}{4096^k}\eq p\(5\l(\f{-1}p\r)-6\l(\f 3p\r)\)\pmod{p^2}.$$
\end{remark}

\noindent{\bf Conjecture II}. {\rm  We have
\begin{align*}\sum_{k=0}^\infty\f{15k+2}{972^k}\bi{2k}k\bi{3k}k
T_k(18,6)=&\f{45\sqrt3}{4\pi},\tag {II1}
\\\sum_{k=0}^\infty\f{91k+12}{10^{3k}}\bi{2k}k\bi{3k}kT_k(10,1)=&\f{75\sqrt3}{2\pi},\tag
{II2}
\\\sum_{k=0}^\infty\f{15k-4}{18^{3k}}\bi{2k}k\bi{3k}kT_k(198,1)=&\f{135\sqrt3}{2\pi},\tag
{II3}
\\\sum_{k=0}^\infty\f{42k-41}{30^{3k}}\bi{2k}k\bi{3k}kT_k(970,1)=&\f{525\sqrt3}{\pi},\tag
{II4}
\\\sum_{k=0}^\infty\f{18k+1}{30^{3k}}\bi{2k}k\bi{3k}kT_k(730,729)=&\f{25\sqrt3}{\pi},\tag {II5}
\\\sum_{k=0}^\infty\f{6930k+559}{102^{3k}}\bi{2k}k\bi{3k}kT_k(102,1)=&\f{1445\sqrt6}{2\pi},\tag
{II6}
\end{align*}
and
\begin{align*}
\\\sum_{k=0}^\infty\f{222105k+15724}{198^{3k}}\bi{2k}k\bi{3k}kT_k(198,1)=&\f{114345\sqrt3}{4\pi},\tag {II7}
\\\sum_{k=0}^\infty\f{390k-3967}{102^{3k}}\bi{2k}k\bi{3k}kT_k(39202,1)=&\f{56355\sqrt3}{\pi},\tag {II8}
\\\sum_{k=0}^\infty\f{210k-7157}{198^{3k}}\bi{2k}k\bi{3k}kT_k(287298,1)=&\f{114345\sqrt{3}}{\pi},\tag {II9}
\\\sum_{k=0}^\infty\f{45k+7}{24^{3k}}\bi{2k}k\bi{3k}kT_k(26,729)=&\f8{3\pi}(3\sqrt3+\sqrt{15}),\tag {II10}
\\\sum_{k=0}^\infty\f{9k+2}{(-5400)^k}\bi{2k}k\bi{3k}kT_k(70,3645)=&\f{15\sqrt3+\sqrt{15}}{6\pi},\tag {II11}
\\\sum_{k=0}^\infty\f{63k+11}{(-13500)^k}\bi{2k}k\bi{3k}kT_k(40,1458)=&\f{25}{12\pi}(3\sqrt3+4\sqrt6).\tag {II12}
\end{align*}}

\begin{remark}\label{Rem5.2}. In view of (44), we may view (II9) as the dual of
(II7) since $198^3/27-198=287298$. The series in (II7) converges
rapidly at a geometric rate with ratio $25/35937$, but the series in
(II9) converges very slow at a geometric rate with ratio
$71825/71874$. (II2), (II9) and (II10)  were motivated by the
following congruences (with $p>3$ a prime):
\begin{align*}&\sum_{k=0}^{p-1}\f{\bi{2k}k\bi{3k}kT_k(10,1)}{10^{3k}}
\\\eq&\begin{cases}4x^2-2p\pmod{p^2}&\t{if}\ p\eq1,3\pmod8\ \&\ p=x^2+2y^2,
\\0\pmod{p^2}&\t{if}\ (\f{-2}p)=-1,\ \t{i.e.},\ p\eq5,7\pmod8,\end{cases}
\end{align*}
\begin{align*}&\sum_{k=0}^{p-1}(91k+12)\f{\bi{2k}k\bi{3k}kT_k(10,1)}{10^{3k}}
\\&\quad\eq\f {3p}2\l(9\l(\f {-3}p\r)-1\r)\pmod{p^2}\ \ (p\not=5);
\end{align*}
\begin{align*}&\sum_{k=0}^{p-1}(210k-7157)\f{\bi{2k}k\bi{3k}kT_k(287298,1)}{198^{3k}}
\\&\quad\eq p\l(35\l(\f {-3}p\r)-7192\r)\pmod{p^2}\ \ \ (p\not=11);
\end{align*}
\begin{align*}&\sum_{k=0}^{p-1}\f{45k+7}{24^{3k}}\bi{2k}k\bi{3k}kT_k(26,729)
\\&\qquad\eq\f p2\l(9\l(\f{-3}p\r)+5\l(\f{-15}p\r)\r)\pmod{p^2}.
\end{align*}
\end{remark}

\noindent{\bf Conjecture III}. {\rm  We have the following formulae:
\begin{align*}\sum_{k=0}^\infty\f{85k+2}{66^{2k}}\bi{4k}{2k}\bi{2k}kT_k(52,1)=&\f{33\sqrt{33}}{\pi},\tag {III1}
\\\sum_{k=0}^\infty\f{28k+5}{(-96^2)^k}\bi{4k}{2k}\bi{2k}kT_k(110,1)=&\f{3\sqrt6}{\pi},\tag {III2}
\\\sum_{k=0}^\infty\f{40k+3}{112^{2k}}\bi{4k}{2k}\bi{2k}kT_k(98,1)=&\f{70\sqrt{21}}{9\pi},\tag {III3}
\\\sum_{k=0}^\infty\f{80k+9}{264^{2k}}\bi{4k}{2k}\bi{2k}kT_k(257,256)=&\f{11\sqrt{66}}{2\pi},\tag {III4}
\\\sum_{k=0}^\infty\f{80k+13}{(-168^2)^k}\bi{4k}{2k}\bi{2k}kT_k(7,4096)=&\f{14\sqrt{210}+21\sqrt{42}}{8\pi},\tag
{III5}
\\\sum_{k=0}^\infty\f{760k+71}{336^{2k}}\bi{4k}{2k}\bi{2k}kT_k(322,1)=&\f{126\sqrt{7}}{\pi},\tag {III6}
\\\sum_{k=0}^\infty\f{10k-1}{336^{2k}}\bi{4k}{2k}\bi{2k}kT_k(1442,1)=&\f{7\sqrt{210}}{4\pi},\tag {III7}
\\\sum_{k=0}^\infty\f{770k+69}{912^{2k}}\bi{4k}{2k}\bi{2k}kT_k(898,1)=&\f{95\sqrt{114}}{4\pi},\tag {III8}
\\\sum_{k=0}^\infty\f{280k-139}{912^{2k}}\bi{4k}{2k}\bi{2k}kT_k(12098,1)=&\f{95\sqrt{399}}{\pi},\tag {III9}
\\\sum_{k=0}^\infty\f{84370k+6011}{10416^{2k}}\bi{4k}{2k}\bi{2k}kT_k(10402,1)=&\f{3689\sqrt{434}}{4\pi},\tag {III10}
\end{align*}
\begin{align*}
\\\sum_{k=0}^\infty\f{8840k-50087}{10416^{2k}}\bi{4k}{2k}\bi{2k}kT_k(1684802,1)=&\f{7378\sqrt{8463}}{\pi},\tag {III11}
\\\sum_{k=0}^\infty\f{11657240k+732103}{39216^{2k}}\bi{4k}{2k}\bi{2k}kT_k(39202,1)=&\f{80883\sqrt{817}}{\pi},\tag {III12}
\\\sum_{k=0}^\infty\f{3080k-58871}{39216^{2k}}\bi{4k}{2k}\bi{2k}kT_k(23990402,1)=&\f{17974\sqrt{2451}}{\pi}.\tag {III13}
\end{align*}}
\begin{remark}\label{Rem5.3} (III12) and (III13) are dual in view of (43). Other
dual pairs include (III6) and (III7), (III8) and (III9), (III10) and
(III11). Below are the corresponding $p$-adic congruences for (III1)
and (III13) (with $p>3$ a prime):
\begin{align*}&\sum_{k=0}^{p-1}(85k+2)\f{\bi{4k}{2k}\bi{2k}kT_k(52,1)}{66^{2k}}
\\\eq& p\l(12\l(\f{-33}p\r)-10\l(\f{33}p\r)\r)\pmod{p^2}\ \ (p\not=11),
\\&\sum_{k=0}^{p-1}(3080k-58871)\f{\bi{4k}{2k}\bi{2k}kT_k(23990402,1)}{39216^{2k}}
\\\eq& p\l(385\l(\f{-2451}p\r)-59256\l(\f{1634}p\r)\r)\pmod{p^2}\ \ (p\not=19,43).
\end{align*}
\end{remark}

\noindent{\bf Conjecture IV}.  We have
\begin{align*}
\sum_{k=0}^\infty\f{26k+5}{(-48^2)^k}\bi{2k}k^2T_{2k}(7,1)=&\f{48}{5\pi},\tag {IV1}
\\\sum_{k=0}^\infty\f{340k+59}{(-480^2)^k}\bi{2k}k^2T_{2k}(62,1)=&\f{120}{\pi},\tag {IV2}
\\\sum_{k=0}^\infty\f{13940k+1559}{(-5760^2)^k}\bi{2k}k^2T_{2k}(322,1)=&\f{4320}{\pi},\tag {IV3}
\end{align*}
\begin{align*}
\\\sum_{k=0}^\infty\f{8k+1}{96^{2k}}\bi{2k}k^2T_{2k}(10,1)=&\f{10\sqrt{2}}{3\pi},\tag {IV4}
\\\sum_{k=0}^\infty\f{10k+1}{240^{2k}}\bi{2k}k^2T_{2k}(38,1)=&\f{15\sqrt6}{4\pi},\tag {IV5}
\\\sum_{k=0}^\infty\f{14280k+899}{39200^{2k}}\bi{2k}k^2T_{2k}(198,1)=&\f{1155\sqrt{6}}{\pi},\tag {IV6}
\\\sum_{k=0}^\infty\f{120k+13}{320^{2k}}\bi{2k}k^2T_{2k}(18,1)=&\f{12\sqrt{15}}{\pi},\tag {IV7}
\\\sum_{k=0}^\infty\f{21k+2}{896^{2k}}\bi{2k}k^2T_{2k}(30,1)=&\f{5\sqrt7}{2\pi},\tag {IV8}
\\\sum_{k=0}^\infty\f{56k+3}{24^{4k}}\bi{2k}k^2T_{2k}(110,1)=&\f{30\sqrt7}{\pi},\tag {IV9}
\\\sum_{k=0}^\infty\f{56k+5}{48^{4k}}\bi{2k}k^2T_{2k}(322,1)=&\f{72\sqrt7}{5\pi},\tag {IV10}
\\\sum_{k=0}^\infty\f{10k+1}{2800^{2k}}\bi{2k}k^2T_{2k}(198,1)=&\f{25\sqrt{14}}{24\pi},\tag {IV11}
\\\sum_{k=0}^\infty\f{195k+14}{10400^{2k}}\bi{2k}k^2T_{2k}(102,1)=&\f{85\sqrt{39}}{12\pi},\tag  {IV12}
\\\sum_{k=0}^\infty\f{3230k+263}{46800^{2k}}\bi{2k}k^2T_{2k}(1298,1)=&\f{675\sqrt{26}}{4\pi},\tag {IV13}
\\\sum_{k=0}^\infty\f{520k-111}{5616^{2k}}\bi{2k}k^2T_{2k}(1298,1)=&\f{1326\sqrt3}{\pi},\tag {IV14}
\\\sum_{k=0}^\infty\f{280k-149}{20400^{2k}}\bi{2k}k^2T_{2k}(4898,1)=&\f{330\sqrt{51}}{\pi},\tag {IV15}
\\\sum_{k=0}^\infty\f{78k-1}{28880^{2k}}\bi{2k}k^2T_{2k}(5778,1)=&\f{741\sqrt{10}}{20\pi},\tag {IV16}
\\\sum_{k=0}^\infty\f{57720k+3967}{439280^{2k}}\bi{2k}k^2T_{2k}(5778,1)=&\f{2890\sqrt{19}}{\pi},\tag {IV17}
\\\sum_{k=0}^\infty\f{1615k-314}{243360^{2k}}\bi{2k}k^2T_{2k}(54758,1)=&\f{1989\sqrt{95}}{4\pi},\tag {IV18}
\end{align*}
and
\begin{align*}\sum_{k=0}^\infty\f{34k+5}{4608^k}\bi{2k}k^2T_{2k}(10,-2)=&\f{12\sqrt6}{\pi},\tag {IV19}
\\\sum_{k=0}^\infty\f{130k+1}{1161216^k}\bi{2k}k^2T_{2k}(238,-14)=&\f{288\sqrt2}{\pi},\tag {IV20}
\\\sum_{k=0}^\infty\f{2380k+299}{(-16629048064)^k}\bi{2k}k^2T_{2k}(9918,-19)=&\f{860\sqrt7}{3\pi}.\tag {IV21}
\end{align*}
\begin{remark}\label{Rem5.4} For (IV6), {\tt Mathematica} indicates that if we set
$$s(n):=\sum_{k=0}^n\f{14280k+899}{39200^{2k}}\bi{2k}{k}^2T_{2k}(198,1)$$
then
$$\bigg|s(15)\times \f{\pi}{1155\sqrt{6}}-1\bigg|<\f1{10^{50}}
\ \ \t{and}\ \ \bigg|s(30)\times \f{\pi}{1155\sqrt{6}}-1\bigg|<\f1{10^{100}}.$$
Here are corresponding $p$-adic congruences of (IV9)-(IV11) and
(IV18) with $p>5$ a prime:
\begin{align*}\sum_{k=0}^{p-1}(56k+3)\f{\bi{2k}k^2T_{2k}(110,1)}{24^{4k}}\eq&\f p4\l(35\l(\f p7\r)-23\r)\pmod{p^2},
\\\sum_{k=0}^{p-1}(56k+5)\f{\bi{2k}k^2T_{2k}(322,1)}{48^{4k}}\eq&\f p{20}\l(147\l(\f p7\r)-47\r)\pmod{p^2},
\end{align*}
\begin{align*}&\sum_{k=0}^{p-1}(10k+1)\f{\bi{2k}k^2T_{2k}(198,1)}{2800^{2k}}
\\\eq&\f p{12}\l(\f 2p\r)\l(13\l(\f p7\r)-1\r)\pmod{p^2}\ \ (p\not=7),
\end{align*}
and
\begin{align*}&\sum_{k=0}^{p-1}(1615k-314)\f{\bi{2k}k^2T_{2k}(54758,1)}{243360^{2k}}
\\\eq&\f p{26}\l(6137\l(\f{p}{95}\r)-14301\r)\pmod{p^2}
\ \ \ (p\not=13).
\end{align*}
For any prime $p>3$, the corresponding $p$-adic congruence of (IV19) is
$$\sum_{k=0}^{p-1}\f{34k+5}{4608^k}\bi{2k}k^2T_{2k}(10,-2)\eq p\l(6\l(\f{-6}p\r)-\l(\f 6p\r)\r)\pmod{p^2}.$$
\end{remark}

\medskip

\noindent{\bf Conjecture V}. {\rm We have the formula
\[\sum_{k=0}^\infty\f{1638k+277}{(-240)^{3k}}\bi{2k}k\bi{3k}kT_{3k}(62,1)=\f{44\sqrt{105}}{\pi}.\tag{V1}\]}

\begin{remark}\label{Rev5.5} (V1) was motivated by Conjecture 4.10; the series
converges at a geometric rate with ratio $-64/125$.
\end{remark}

We conjecture that (IV1)-(IV18) have exhausted all identities of the
form
$$\sum_{k=0}^\infty(a+dk)\f{\bi{2k}k^2T_{2k}(b,1)}{m^k}=\f C{\pi}$$
with $a,d,m\in\Z$, $b\in\{1,3,4,\ldots\}$, $d>0$, and $C^2$ rational. This comes from our following hypothesis
motivated by (\ref{2.5}) in the case $h=2$ and the author's philosophy
about series for $1/\pi$ stated in \cite{Su7}. We have applied the
hypothesis to seek for series for $1/\pi$ of type IV and checked all
those $b=1,\ldots,10^6$ via computer.
\medskip

\noindent{\bf Hypothesis 5.1} {\rm (i) Suppose that
$$\sum_{k=0}^\infty\f{a+dk}{m^k}\bi{2k}k^2T_{2k}(b,1)=\f C{\pi}$$
with $a,d,m\in\Z$, $b\in\Z^+$ and $C^2\in\Q\sm\{0\}$. Then
$\sqrt{|m|}$ is an integer dividing $16(b^2-4)$. Also, $b=7$ or
$b\eq2\pmod4$.

 (ii) Let $\ve\in\{\pm1\}$, $b,m\in\Z^+$ and $m\mid 16(b^2-4)$.
Then, there are $a,d\in\Z$ such that
$$\sum_{k=0}^\infty\f{a+dk}{(\ve m^2)^k}\bi{2k}k^2T_{2k}(b,1)=\f C{\pi}$$
for some $C\not=0$ with $C^2$ rational, if and only if $m>4(b+2)$
and
$$\sum_{k=0}^{p-1}\f{\bi{2k}k^2T_{2k}(b,1)}{(\ve m^2)^k}
\eq\l(\f{\ve(b^2-4)}p\r)\sum_{k=0}^{p-1}\f{\bi{2k}k^2T_{2k}(b,1)}{(\ve
\bar m^2)^k}\pmod{p^2}$$ for all odd primes $p\nmid b^2-4$, where
$\bar m=16(b^2-4)/m$.}

\medskip

\noindent{\bf Conjecture VI}. {\rm We have the following formulae:
\begin{align*}
\sum_{k=0}^\infty\f{66k+17}{(2^{11}3^3)^{k}}T_k^3(10,11^2)=&\f{540\sqrt2}{11\pi},\tag {VI1}
\\\sum_{k=0}^\infty\f{126k+31}{(-80)^{3k}}T_k^3(22,21^2)=&\f{880\sqrt5}{21\pi},\tag {VI2}
\\\sum_{k=0}^\infty\f{3990k+1147}{(-288)^{3k}}T_k^3(62,95^2)=&\f{432}{95\pi}(195\sqrt{14}+94\sqrt2).\tag {VI3}
\end{align*}}

\begin{remark}\label{Rem5.5} The series (VI1)-(VI3) converge at geometric rates with ratios
$$\f{16}{27},\ -\f{64}{125},\ -\f{343}{512}$$
respectively. The author would like to offer \$300 as the prize for the person (not joint authors)
who can provide first rigorous proofs of all the three identities (VI1)-(VI3).
(VI1) and (VI3) were motivated by the author's following conjectural congruences for any prime $p>3$:
$$\sum_{k=0}^{p-1}\f{T_k^3(10,11^2)}{(2^{11}3^3)^k}\eq\begin{cases}(\f{2}p)(4x^2-2p)\,(\mo\ p^2)&\t{if}\ p=x^2+3y^2,
\\0\ (\mo\ p^2)&\t{if}\ p\eq2\ (\mo\ 3),\end{cases}$$
$$\sum_{k=0}^{p-1}\f{66k+17}{(2^{11}3^3)^k}T_k^3(10,11^2)
\eq\f p{11}\l(195\l(\f{-2}p\r)-8\l(\f{-6}p\r)\r)\pmod{p^2};$$
$$\sum_{k=0}^{p-1}\f{T_k^3(62,95^2)}{(-288)^{3k}}\eq\begin{cases}(\f{-2}p)(4x^2-2p)\,(\mo\ p^2)&\t{if}\ p=x^2+7y^2,
\\0\ (\mo\ p^2)&\t{if}\ (\f p7)=-1,\end{cases}$$
\begin{align*}&\sum_{k=0}^{p-1}\f{3990k+1147}{(-288)^{3k}}T_k^3(62,95^2)
\\\eq&\f p{19}\l(17563\l(\f{-14}p\r)+4230\l(\f{-2}p\r)\r)\pmod{p^2}.
\end{align*}
\end{remark}

\medskip

\noindent{\bf Conjecture VII}. {\rm We have the following formulae:
\begin{align*}
\sum_{k=0}^\infty\f{221k+28}{450^{k}}\bi{2k}kT_k^2(6,2)=&\f{2700}{7\pi},\tag {VII1}
\\\sum_{k=0}^\infty\f{24k+5}{28^{2k}}\bi{2k}kT_k^2(4,9)=&\f{49}{9\pi}(\sqrt3+\sqrt6),\tag {VII2}
\\\sum_{k=0}^\infty\f{560k+71}{22^{2k}}\bi{2k}kT_k^2(5,1)=&\f{605\sqrt7}{3\pi},\tag {VII3}
\\\sum_{k=0}^\infty\f{3696k+445}{46^{2k}}\bi{2k}kT_k^2(7,1)=&\f{1587\sqrt7}{2\pi},\tag{VII4}
\\\sum_{k=0}^\infty\f{56k+19}{(-108)^k}\bi{2k}kT_k^2(3,-3)=&\f{9\sqrt7}{\pi},\tag {VII5}
\\\sum_{k=0}^\infty\f{450296k+53323}{(-5177196)^k}\bi{2k}kT_k^2(171,-171)=&\f{113535\sqrt7}{2\pi},\tag {VII6}
\\\sum_{k=0}^\infty\f{2800512k+435257}{434^{2k}}\bi{2k}kT_k^2(73,576)=&\f{10406669}{2\sqrt6\,\pi}.\tag {VII7}
\end{align*}}

\begin{remark}\label{Rem5.6} The series (VII1)-(VII7) converge at geometric rates with ratios
$$\f{88+48\sqrt2}{225},\ \f{25}{49},\ \f{49}{121},\ \f{81}{529},\ -\f 79,\ -\f{175}{7569},\ \f{14641}{47089}.$$
respectively. The author found (VII2) and (VII3) in light of Conjecture 4.18. Similarly, (VII6)-(VII7) were motivated by Conjectures 4.19-4.20.
\end{remark}

Concerning the new identities in Conjectures I-VII, actually we
first discovered congruences without linear parts related to binary
quadratic forms (like many congruences in Section 4), then found
corresponding $p$-adic congruences with linear parts, and finally
figured out the series for $1/\pi$.

\Ack. The work was supported by the National Natural Science Foundation (grant 11171140) of China, and
the initial version of this paper was posted to {\tt arXiv} in Jan. 2011 as a preprint with the ID {\tt arXiv:1101.0600}.
The preprint version of this paper available from {\tt arXiv} has stimulated some others to work on our conjectural series for $1/\pi$
of types I-V in Section 5.

\end{document}